\documentclass[]{amsart}
\usepackage[
  top=0.95in,
  bottom=0.85in,
  left=0.90in,
  right=0.90in,
  bindingoffset=0.25in,
  heightrounded,
]{geometry}
\usepackage{babel}
\usepackage[utf8]{inputenc}
\usepackage{amsmath}
\usepackage{amsthm}
\usepackage{amsfonts}
\usepackage{amssymb}
\usepackage{enumerate}
\usepackage{stackrel}
\usepackage{mathrsfs}
\usepackage{graphicx}
\usepackage{tikz-cd}
\usepackage[hidelinks]{hyperref}
\usepackage{fancyhdr}
\usepackage{cleveref}
\usetikzlibrary{babel}
\author{Iván Chércoles Cuesta}
\title{Spaces of left-preorders on free products}

\hypersetup{
     colorlinks=true,
     linkcolor=blue,
     citecolor=red,
     }

\providecommand{\df}[3]{#1:#2 \longrightarrow #3}
\providecommand{\natp}{\mathbb{N}_{>0}}
\providecommand{\nat}{\mathbb{N}}
\providecommand{\ent}{\mathbb{Z}}
\providecommand{\real}{\mathbb{R}}
\providecommand{\co}{\circ}
\providecommand{\sii}{\quad \Longleftrightarrow \quad}

\providecommand{\cont}{\subseteq}
\providecommand{\inv}{^{-1}}
\providecommand{\tq}{\;/\;}

\providecommand{\gen}[1]{\langle #1\rangle}
\providecommand{\semgen}[1]{\langle #1\rangle^{+}}
\providecommand{\rela}{\mathcal{R}}

\providecommand{\homeomr}{\mathrm{Homeo}_{+}(\real)}
\providecommand{\lleq}{\preceq}
\providecommand{\lle}{\prec}
\providecommand{\lleqb}{\lleq}

\providecommand{\settminus}{-}

\providecommand{\stab}{\mathrm{Stab}}

\providecommand{\pord}{\mathcal{PO}}

\providecommand{\gr}{G}

\providecommand{\grb}{H}

\providecommand{\ga}{g}
\providecommand{\gb}{h}

\providecommand{\gy}{y}
\providecommand{\gx}{x}

\providecommand{\ce}{C}

\providecommand{\grf}{\Gamma}
\providecommand{\grfb}{\Omega}
\providecommand{\grfv}{\vertices \grf}
\providecommand{\grfe}{\edges \grf}
\providecommand{\grfbv}{\vertices \grfb}
\providecommand{\grfbe}{\edges \grfb}
\providecommand{\grfbv}{\vertices \grfb}
\providecommand{\grfbe}{\edges \grfb}
\providecommand{\pa}{\mathtt{p}}

\providecommand{\vea}{\mathtt{v}}
\providecommand{\veb}{\mathtt{u}}

\providecommand{\equis}{\mathcal{X}}
\providecommand{\equisg}{\equis_{\gr}}
\providecommand{\equish}{\equis_{\grb}}

\providecommand{\origen}{\mathtt{o}}
\providecommand{\final}{\mathtt{t}}
\providecommand{\eda}{\mathtt{e}}
\providecommand{\edb}{\mathtt{f}}
\providecommand{\arbol}{\mathtt{T}}
\providecommand{\edges}{\mathtt{E}}
\providecommand{\vertices}{\mathtt{V}}
\providecommand{\veagr}{\vea_{\gr}}
\providecommand{\veacero}{\vea_0}
\providecommand{\veagrb}{\vea_{\grb}}
\providecommand{\edagr}{\eda_{\gr}}
\providecommand{\edagrb}{\eda_{\grb}}
\providecommand{\edagr}{\eda_{\gr}}
\providecommand{\camaso}{\mathfrak{p}}
\providecommand{\vebase}{\veb_*}
\providecommand{\fino}{\Omega}
\providecommand{\fini}{I}
\providecommand{\finj}{J}
\providecommand{\fink}{K}
\providecommand{\indi}{i}
\providecommand{\indj}{j}
\providecommand{\indk}{k}
\providecommand{\arbmax}{\arbol_0}
\providecommand{\aristasfuera}{\mathtt{Y}}
\providecommand{\generadores}{S}

\providecommand{\entorno}{\mathbb{B}}

\providecommand{\cete}{c}
\providecommand{\rank}{\mathit{rank}}

\providecommand{\defrelativeorder}[2]{left-preorder on $#1$ relative to $#2$}

\providecommand{\conjfiny}{\mathscr{T}}
\providecommand{\elefiny}{t}
\providecommand{\conjfinf}{\mathscr{R}}
\providecommand{\elefinf}{r}
\providecommand{\puntox}{\mathsf{x}}
\providecommand{\puntoy}{\mathsf{y}}

\providecommand{\gfa}{a}
\providecommand{\gfb}{b}
\providecommand{\grfaca}{A}
\providecommand{\grfacb}{B}

\providecommand{\bola}{\entorno}

\newtheorem{ej}{Example}[section]
\newtheorem{teo}[ej]{Theorem}
\newtheorem{prop}[ej]{Proposition}
\newtheorem{lema}[ej]{Lemma}
\newtheorem{coro}[ej]{Corollary}
\theoremstyle{definition}
\newtheorem{defi}[ej]{Definition}
\newtheorem{obs}[ej]{Remark}
\newtheorem*{nota}{Notation}

\usepackage{blindtext}

\begin{document}

\address{Instituto de Ciencias Matemáticas (CSIC-UAM-UC3M-UCM), Consejo Superior de Investigaciones Científicas, C/ Nicolás Cabrera, 13–15, Campus de Cantoblanco UAM, 28049 Madrid, Spain.}
\email{ivan.chercoles@gmail.com}
\date{\today}

\let\thefootnote\relax
\footnotetext{The author is supported by the grant CEX2019-000904-S-20-1 funded by: MCIN/AEI/ 10.13039/501100011033. The author is also partially supported by the grant PID2021-126254NB-I00 of the Ministry of Science and Innovation of Spain.}

\maketitle

\begin{abstract}
Using dynamical techniques we show that there are no isolated elements on the space of left-preorders on a free product of two groups. As a consequence, when the groups are finitely generated, this space is either empty or a Cantor set. For any subgroup of a free product having finite Kurosh rank, we develop similar results for the subspace consisting on the set of left-preorders relative to that subgroup. We provide a generating set for a certain intersection of subgroups of a free product.
\end{abstract}


\section{Introduction}

In this paper we study the space of left-preorders on a free product. We start by recalling the definition of left-preorder on a group.

\begin{defi}
Let $\gr$ be a group. A \textit{left-preorder} $\lleq$ on $\gr$ is a total transitive reflexive binary relation that is left-invariant, meaning that $$\ga_1\lleq \ga_2\mbox{ implies }\ga \ga_1\lleq \ga \ga_2\mbox{  for all }\ga,\ga_1,\ga_2\in \gr.$$ The \textit{trivial left-preorder} is defined by $\ga_1\lleq \ga_2$ for all $\ga_1,\ga_2\in \gr$. We define $\pord(\gr)$ the \textit{set of non-trivial left-preorders} on $\gr$. For convenience, we are going to assume that the trivial left-preorder is not a left-preorder.
\end{defi}

Any left-preorder $\lleq$ on $\gr$ give us a unique decomposition $\gr=P\sqcup P\inv \sqcup \ce$ where $P$ is the positive cone of $\lleq$ and $\ce$ is a proper subgroup. Then, we can associate $\lleq$ with the map 
\begin{equation*}
\begin{array}{r@{\hspace{0pt}}c@{\hspace{0pt}}
c@{\hspace{4pt}}l}
\phi_\lleq:&\gr &\longrightarrow& \{\pm1,0\} \\

&\ga\ \ &\mapsto&\begin{cases}1 \quad \textrm{if} \quad \ga\in P \\
 -1 \quad \textrm{if}\quad \ga\in P\inv \\
 0 \quad \textrm{otherwise} .\end{cases}
\end{array}
\end{equation*} 
This give us an injective map from $\pord(\gr)$ to $\{\pm 1,0\}^\gr$ by sending each left-preorder $\lleq$ to the map $\phi_{\lleq}$. Hence, we can identify $\pord(\gr)$ with a subset of the set $\{\pm 1,0\}^\gr$. We consider on $\{\pm 1,0\}$ the discrete topology and we consider on $\{\pm 1,0\}^\gr$ the product topology.

\begin{defi}
The \textit{space of left-preorders} of $\gr$ is defined to be the set $\pord(\gr)$ endowed with the subspace topology of $\{\pm 1,0\}^\gr$ with the product topology.
\end{defi}

The space of left-preorders of a finitely generated group is always metrizable, compact and totally disconnected.\\

The \textit{space of left-orders} is the subspace of the space of left-preorders given by the set of antisymmetric left-preorders. McClearly proved in \cite{free-mcclearly} that there are no isolated left-orders on any free group of finite rank greater than 1 so, in particular, its space of left-orders is a Cantor set. Rivas in \cite{freeproducts-rivas} generalized this fact for free products of two non-trivial finitely generated left-orderable groups. It is also remarkable that in \cite{biorderablefreeproduct} Muliarchyk showed that the space of bi-orders (the subspace of the space of left-orders given by the ones that are also right-invariant) of a free product of two non-trivial finitely generated biorderable groups is a Cantor set. We prove the following theorem. 

\begin{teo}[\Cref{teorema aislados general} \& \Cref{corolario espacio de preordenes es cantor producto libre}]\label{teorema aislados general introduccion}
Let $\gr$ and $\grb$ be two non-trivial left-orderable groups. Then $\pord(\gr*\grb)$ has no isolated elements. In particular, if $\gr$ and $\grb$ are finitely generated, then $\pord(\gr*\grb)$ is a Cantor set.
\end{teo}

As a corollary, we deduce that the space of left-preorders of a finitely generated non-abelian free group is a Cantor set.\\

The space of left-preorders has been studied by Decaup and Rond in \cite{valuated-decauprond} for the abelian case. This space has some applications for other subjects of mathematics, for instance, in \cite{locallymoving} it is shown that a certain natural quotient on a product of the space of left-preorders it is homeomorphic to a space of normalized harmonic actions, in \cite{alggeo1} it plays a central role to obtain completions of real fans, and in  \cite{alggeo2} it is considered to study toric varieties. In \cite{morris-antolinrivas} an analogue of Morris theorem \cite{amenable-morris} is proven for left-preoders using the compactness of the space of left-preorders. The space of conradian left-preorders, which is a subspace of the space of left-orders, it is studied in \cite{chercoconrad}.\\

We now briefly describe the strategy used on the proof of the theorem. Before we need a definition. Recall that every left-preorder $\lleq$ on a group $\gr$ give us a unique decomposition $\gr=P\sqcup P\inv \sqcup \ce$ where $\ce$ is a proper subgroup, and we will say that $\lleq$ is a left-preorder on $\gr$ \textit{relative to} $\ce$.

\begin{defi}
Let $\gr$ be a group and $\ce$ a proper subgroup of $\gr$. We say that $\ce$ is \textit{left-relatively convex} on $\gr$ if there exists some left-preorder on $\gr$ relative to $\ce$. We define $\pord_\ce(\gr)$ to be the \textit{set of left-preorders on} $\gr$ \textit{relative to} $\ce$. The \textit{space of left-preorders on} $\gr$ \textit{relative to} $\ce$ is defined to be $\pord_\ce(\gr)$ endowed with the topology given as a subspace of $\pord(\gr)$.
\end{defi}

The space $\pord(\gr)$ can be decomposed as a union of $\pord_\ce(\gr)$ where $\ce$ ranges over all left-relatively convex proper subgroups on $\gr$. In order to derive the theorem, we show the following theorem about the structure of $\pord_\ce(\gr)$. 

\begin{teo}[\Cref{teorema aislados} \& \Cref{espacio de ordenes es un conjunto de Cantor}]\label{teorema relativo aislados introduccion}
Let $\gr$ and $\grb$ be two non-trivial groups. Let $\ce$ be a left-relatively convex proper subgroup. If $\ce$ is finitely generated, then $\pord_\ce(\gr*\grb)$ has no isolated elements.
\end{teo}

The proof uses the same strategy done by Rivas in \cite{freeproducts-rivas}. We recover Rivas' result about the space of left-orders when taking $\ce=\{1\}$ in the prevous theorem. Actually, we prove a more general result: in last theorem we can change the condition of $\ce$ being finitely generated by having finite Kurosh rank on the free product $\gr*\grb$, which is a weaker condition. The \textit{Kurosh rank} of a subgroup $\ce$ of $\gr*\grb$ is defined as the number of free factors in the decomposition of $\ce$ according to Kurosh's subgroup theorem. Equivalently, the Kurosh rank of $\ce$ is the sum of the number of vertices of $\ce\backslash\arbol$ with associated non-trivial stabilizer plus the number of edges of $\ce\backslash\arbol$ outside of a maximal subtree, for $\arbol$ the Bass-serre universal covering tree relative to the free product $\gr*\grb$.\\

The proof of the theorem needs a result providing a generating set for an intersection of subgroups of a free product. This fact, as far as the author knows, has not been considered in the literature, so we state it as a result that could be interesting by itself.

\begin{teo}[\Cref{bass serre para interseccion}]
Let $\ce$ be a proper subgroup of $\gr*\grb$. If $\ce$ has finite Kurosh rank, then there exists $\fino,\{\zeta_\indi\}_{\indi\in\fini},\{\alpha_\indj\}_{\indj\in\finj},\{\beta_\indk\}_{\indk\in\fink}\subseteq \gr*\grb$ finite sets satisfying the following property:\\

For every $\gr_0$ subgroup of $\gr$ and every $\grb_0$ subgroup of $\grb$ satisfying $\fino\subseteq \gen{\gr_0,\grb_0}$ there exists $\gr_\indj$ subgroup of $\gr_0$ for each $\indj\in\finj$ and there exists $\grb_\indk$ subgroup of $\grb_0$ for each $\indk\in\fink$ such that 
\begin{equation*}
\ce\cap\gen{\gr_0,\grb_0}=\gen{\{\zeta_\indi \tq \indi\in\fini\}\cup \left(\bigcup_{\indj\in\finj}\alpha_\indj\gr_\indj\alpha_\indj\inv \right)\cup \left(\bigcup_{\indk\in\fink}\beta_\indk\grb_\indk\beta_\indk\inv \right)}.
\end{equation*}
\end{teo}

It is remarkable that, for such $\ce$ and for such $\gr_0$ and $\grb_0$ satisfying $\fino\subseteq \gen{\gr_0,\grb_0}$, we have that the Kurosh rank of $\ce\cap\gen{\gr_0,\grb_0}$ is bounded by the Kurosh rank of $\ce$.\\

We also mention a relevant consequence of the theorem about $\pord_\ce(\gr*\grb)$. Recall that the positive cone of a left-preorder is always a semigroup. The fact that $\pord_\ce(\gr*\grb)$ has no isolated elements implies that the positive cone of any left-preorder relative to $\ce$ is not finitely generated as a semigroup. In particular, the mentioned theorem implies the following. 

\begin{coro}[\Cref{cono positivo no finitamente generado como corolario de aislado}]\label{corolario parcial cono positivo finitamente generado introduccion}
Let $\gr$ and $\grb$ be non-trivial finitely generated groups. Let $\ce$ be a left-relatively convex proper subgroup. If $\ce$ is finitely generated, then no left-preorder on $\gr*\grb$ relative to $\ce$ has a positive cone $P$ that is finitely generated as a subsemigroup.
\end{coro}

On the previous corollary, one can change the condition of $\ce$ being finitely generated by the condition of $\ce$ having finite Kurosh rank.

\section{Left-preorders and positive cones}

We start by introducing left-preorders. Left-preorders are also considered in other works, for example in \cite{chercoconrad} with the same 
nomenclature, in \cite{locallymoving, valuated-decauprond} with the name of \textit{left-invariant preorder}, and in \cite{morris-antolinrivas} with the name of \textit{relative order}. We will assume that every left-preorder is not trivial. Our notation of left-preorders is the one in \cite{chercoconrad}, we detail it here in order to be clear with the exposition.\\

A \textit{preorder} $\lleq$ on a set $X$ is a transitive reflexive binary relation that is total (for every $x_1,x_2\in X$ we have $x_1\lleq x_2$ or $x_2\lleq x_1$). We define $[x]_{\lleq}$ as the set of $y\in X$ such that $x\lleq y$ and $y\lleq x$. A preorder $\lleq$ on $X$ is \textit{trivial} if $x_1\lleq x_2$ for all $x_1,x_2\in X$. An order on $X$ is a preorder that is antisymmetric.\\

Given a group $\gr$ acting on a set $X$, a relation $\rela$ on $X$ is \textit{invariant under the action} of $\gr$ if $x_1\rela x_2$ implies $\ga x_1\rela \ga x_2$ for every $x_1,x_2\in X$ and all $\ga\in \gr$. A left-preorder $\lleq$ on a group $\gr$ is a \textit{left-preorder} on $\gr$ that is non-trivial and that is invariant under the action of $\gr$ on itself by left multiplication $$\ga_1\lleq \ga_2\mbox{ implies }\ga \ga_1\lleq \ga \ga_2\mbox{  for all }\ga,\ga_1,\ga_2\in \gr.$$ A \textit{left-order} on $\gr$ is an antisymmetric left-preorder on $\gr$.\\

We state a lemma that allows us to identify left-preorders on $\gr$ with total left-invariant orders on sets of cosets. We do not write the proof as it is routinary. Notice that a left-preorder $\lleq$ on $\gr$ is a left-order on $\gr$ if and only if $[1]_{\lleq}=\{1\}$. In next lemma, when $\ce$ is normal, the total order invariant under left multiplication on $\gr/\ce$ is a left-order on $\gr/\ce$.

\begin{lema}\label{lema preorden como orden en cocientes}
Let $\gr$ be a group. Given $\lleq$ a left-preorder on $\gr$, then $[1]_{\lleq}$ is a proper subgroup of $\gr$ and the relation $\lleqb$ on $\gr/[1]_{\lleq}$ defined by $$\ga[1]_{\lleq}\lleqb \gb[1]_{\lleq}\; \Leftrightarrow\; \ga\lleq \gb $$ is a total order invariant under left multiplication. Reciprocally, given $\ce$ a proper subgroup of $\gr$ and $\lleqb$ a total order on $\gr/\ce$ invariant under left multiplication, the relation $\lleq$ on $\gr$ defined by $$\ga\lleq \gb \; \Leftrightarrow\; \ga \ce\lleqb \gb \ce$$ is a left-preorder on $\gr$ satisfying $\ce=[1]_{\lleq}$.
\end{lema} 

\begin{nota}
We identify every left-preorder on $\gr$ with a total left-invariant order on a set of cosets over a subgroup of $\gr$. If $\lleq$ is a left-preorder on $\gr$, $\ga [1]_{\lleq}\lle \gb [1]_{\lleq}$ will mean that $\ga [1]_{\lleq}\lleq \gb [1]_{\lleq}$ and $\ga [1]_{\lleq}\neq \gb [1]_{\lleq}$.
\end{nota}

A left-preorder on a group $\gr$ \textit{relative to} a proper subgroup $\ce$ is a left preorder $\lleq$ satisfying that $\ce=[1]_{\lleqb}$. We denote by $\pord_\ce(\gr)$\textit{ the set of left-preorders on} $\gr$ \textit{relative to} $\ce$. We say that $\ce$ is \textit{left-relatively convex} on $\gr$ if there is a left-preorder on $\gr$ relative to $\ce$.

The notion of left-relatively convex is considered in \cite{relativelyconvex-antolinwarrenzoran}. In that reference it is described some examples of left-relatively convex subgroups, for instance maximal cyclic subgroups of non-abelian free groups or maximal $n$-generated abelian subgroups of right-angled Artin groups. From Theorem 17 in \cite{relativelyconvex-antolinwarrenzoran} we have the following fact that we will need later.

\begin{prop}\label{ordenar con factor convexo}
Let $(\gr_i)_{i\in I}$ be a family of non-trivial left-orderable groups. Take $k\in I$. Then, $\gr_{k}$ is a left-relatively convex subgroup on $\ast_{i\in I} \gr_i$.
\end{prop}

We now consider positive cones. We do not prove next proposition as it is routinary with the ideas of Definition 4 on \cite{relativelyconvex-antolinwarrenzoran}.

\begin{nota}
Let $A,B$ be two sets. We denote $A\sqcup B$ the union $A\cup B$ if we also have that $A$ and $B$ are disjoint.
\end{nota}

\begin{prop}\label{orden relativo como semigrupo}
Given $\lleq$ a \defrelativeorder{\gr}{\ce}, the set $P=\{\ga\in \gr\tq \ce \lle \ga \ce \}$ satisfies both
\begin{enumerate}
\item[a)] $P$ is a subsemigroup of $\gr$.
\item[b)] $\gr=P\sqcup P\inv\sqcup \ce $.
\end{enumerate}
Such $P$ is called the \textit{positive cone} of $\lleq$. Reciprocally, if $\ce $ is a proper subgroup of $\gr$ and $P\cont \gr$ is a subset in such a way the last two conditions are satisfied, defining $$\ga\ce \lleq \gb\ce  \Leftrightarrow \ga\inv \gb\in P\cup \ce $$ give us a left-preorder on $\gr$ relative to $\ce $. 
\end{prop}

\begin{obs}\label{positivo por neutro es positivo}
Given $\lleq$ a \defrelativeorder{\gr}{\ce}, the positive cone $P$ satisfies $P\ce\subseteq P$ and $\ce P\subseteq P$.
\end{obs}

We now describe left-preorders from a dynamical point of view. For a complete reference of these dynamical techniques see \cite{god}. From now until the end of this section, we will follow the ideas and the notation from Section 1 of \cite{freeproducts-rivas} but adapted from left-orders to left-preorders.\\

We will denote by $\homeomr$ the set of orientation preserving homeomorphisms of the real line. Every orientation preserving homeomorphism of the real line preserves the usual order on $\real$. In particular, $\homeomr$ is a subgroup of the group of order preserving bijections of the real line endowed with the usual order.

\begin{nota}
Let $\gr$ be a group. Each homomorphism $\df{\rho}{\gr}{\homeomr}$ induces an action of $\gr$ on $\real$ by orientation preserving homeomorphisms and vice versa. Given an homomorphism $\df{\rho}{\gr}{\homeomr}$ and $p\in\real$, we will denote by $\stab_\rho(p)$ the stabilizer of $p$ with respect to the action induced by the homomorphism $\rho$. In other words, $$\stab_\rho(p)=\{\ga \in \gr \tq \rho(\ga)(p)=p\}.$$
\end{nota}

The following proposition can be found in \cite{morris-antolinrivas} as Proposition 2.2, so we do not prove it.

\begin{prop}\label{lema realizacion dinamica orden}
Let $\gr$ be a group and $\ce$ a proper subgroup.
If there is a homomorphism $\df{\rho}{\gr}{\homeomr}$ and $p\in \real$ such that $\stab_\rho(p)=\ce$, then defining 
$$\ga_1 \ce\lleq_{\rho} \ga_2 \ce \Longleftrightarrow \rho(\ga_1)(p)\leq\rho(\ga_2)(p)$$
give us a left-preorder $\lleq_{\rho}$ on $\gr$ relative to $\ce$. Reciprocally, if $\gr/\ce$ is countable and $\lleq$ is a left-preorder on $\gr$ relative to $\ce$, then there is a homomorphism $\df{\rho_{\lleq}}{\gr}{\homeomr}$ and an element $p\in \real$ such that $\stab_{\rho_{\lleq}}(p)=\ce$ and 
$$\ga_1 \ce\lleq \ga_2 \ce \Longleftrightarrow \rho_{\lleq}(\ga_1)(p)\leq\rho_{\lleq}(\ga_2)(p).$$
We say that $\rho_{\lleq}$ is a \textit{dynamical realization-like homomorphism} for $\lle$ with reference point $p$. Moreover, we can build $\rho_{\lleq}$ in such a way $p=0$.
\end{prop}

\begin{obs}
Let $\df{D}{\gr}{\homeomr}$ be a dynamical realization-like homomorphism for $\lleq$ with reference point $p$ and let $\df{\phi}{\real}{\real}$ be an order preserving homeomorphism. We can define the conjugation $\df{D_\phi}{\gr}{\homeomr}$ given by $D_\phi(g)=\phi \co D(g)\co \phi\inv$ for each $g\in \gr$. Then, $D_\phi$ is a dynamical realization-like homomorphism for $\lleq$ with reference point $\phi(p)$.
\end{obs}

We now introduce some notation and some technical conclusions that we will need to prove \Cref{teorema aislados}. Given $\equis_{\gr}$ a finite subset of  group $\gr$, we denote by $\bola_{\equis_{\gr}}(n)$ the ball of length $n$ with respect to $\equis_{\gr}$, that is, all the elements that result just by multiplying $n$ or less elements in $\equis_{\gr}$.

\begin{defi}
Let $\gr$ be a group and $\equis_{\gr}$ a finite subset such that $\equis_{\gr}\inv=\equis_{\gr}$. Let $\lleq$ be a left-preorder on $\gr$ relative to $\ce$. We say that an element $\ga\in \gr$ is $\lleq$\textit{-maximal element in the ball} $\bola_{\equis_{\gr}}(n)$ if $$\ga\in \bola_{\equis_{\gr}}(n) \;\mbox{ and }\;\ga \ce=\max_{\lleq}\{w\ce\tq w\in \bola_{\equis_{\gr}}(n)\}.$$ We say that an element $\ga\in \gr$ is $\lleq$\textit{-minimal element in the ball} $\bola_{\equis_{\gr}}(n)$ if $$\ga\in \bola_{\equis_{\gr}}(n)\; \mbox{ and }\;\ga \ce=\min_{\lleq}\{w\ce\tq w\in \bola_{\equis_{\gr}}(n)\}.$$ 
\end{defi}

Assuming the notation of the previous definition, we have that there always exists a  $\lleq$-maximal element in the ball $\bola_{\equis_{\gr}}(n)$ and a $\lleq$-minimal element in the ball $\bola_{\equis_{\gr}}(n)$. This is a direct consequence of the fact that $\bola_{\equis_{\gr}}(n)$ is finite since $\equis_{\gr}$ is finite. Also notice that these elements might not be unique. However, given $D$ a dynamical realization-like homomorphism for $\lleq$ with reference point $p$, the value $D(\lambda_n^+)(p)$ do not depend on the chosen $\lambda_n^+$ $\lleq$-maximal element in the ball $\bola_{\equis_{\gr}}(n)$, since $\stab_D(p)=\ce$ (and similarly for the $\lleq$-minimal case). 

\begin{lema}\label{lema comparando lambdas}
Let $\gr$ be a group and $\equis_{\gr}$ a finite subset such that $\equis_{\gr}\inv=\equis_{\gr}$. Let $\lleq$ be a left-preorder on $\gr$ relative to $\ce$. Assume that $\ce\cap \gen{\equis_{\gr}}\neq \gen{\equis_{\gr}}$. For each $n\in\nat$, consider $\lambda_n^+$ (resp. $\lambda_n^-$) a $\lleq$-maximal (resp. $\lleq$-minimal) element in the ball $\bola_{\equis_{\gr}}(n)$. Then, we have $\lambda_n^+ \ce\lle \lambda_{n+1}^+\ce$ and $\lambda_{n+1}^-\ce\lle \lambda_n^-\ce$ for all $n\in\nat$. In particular, $\lambda_n^+\ce\neq\lambda_k^+\ce$ and $\lambda_n^-\ce\neq\lambda_k^-\ce$ for all $n,k\in\nat$ such that $n\neq k$.
\begin{proof}
We prove the first equation $\lambda_n^+ \ce\lle \lambda_{n+1}^+\ce$ for all $n\in\nat$, as the proof of $\lambda_{n+1}^-\ce\lle \lambda_n^-\ce$ for all $n\in\nat$ is analogous. Fix $n\in\nat$. Since $\ce\cap \gen{\equis_{\gr}}\neq \gen{\equis_{\gr}}$ and since $\equis_{\gr}\inv=\equis_{\gr}$, we can take $\gy \in \equis_{\gr}$ such that $\ce\lle\gy\ce$. In particular, $\lambda_n^+\ce\lle\lambda_n^+\gy\ce$. Then, since $\lambda_n^+\gy\in \bola_{\equis_{\gr}}(n+1)$, we obtain $\lambda_n^+\ce\lle\lambda_n^+\gy\ce\lleq\lambda_{n+1}^+\ce$.
\end{proof}
\end{lema}

\begin{lema}\label{lema realizaciones}
Let $\gr$ be a countable group and $\equis_{\gr}$ a finite subset such that $\equis_{\gr}\inv=\equis_{\gr}$. Consider $\ce$ a left-relatively convex subgroup on $\gr$ and $\lleq$ a left-preorder on $\gr$ relative to $\ce$. Let $\df{D}{\gr}{\homeomr}$ be a dynamical realization-like homomorphism for $\lleq$ with reference point $p$ and $n\in\nat$. For each $n\in\nat$, consider $\lambda_n^+$ a $\lleq$-maximal element in the ball $\bola_{\equis_{\gr}}(n)$ and consider $\lambda_n^-$ a $\lleq$-minimal element in the ball $\bola_{\equis_{\gr}}(n)$. Then,
\begin{enumerate}
\item For each $\ga\in \bola_{\equis_{\gr}}(n)$, we have that $D(\ga)(p)$ lies in the interval $[D(\lambda_n^-)(p),D(\lambda_n^+)(p)]$.
\item Suppose that $\df{\tilde{D}}{\gr}{\homeomr}$ is a homomorphism. If for all $\ga\in \equis_{\gr}$ we have that $D(\ga)$ and $\tilde{D}(\ga)$ are equal as maps when restricted to $[D(\lambda_n^-)(p),D(\lambda_n^+)(p)]$, then $D(\ga)(p)=\tilde{D}(\ga)(p)$ for all $\ga\in \bola_{\equis_{\gr}}(n+1)$.
\end{enumerate}
\begin{proof}
We prove each statement. We start by proving the first one. Take $\ga\in \bola_{\equis_{\gr}}(n)$. By definition, we have $\lambda_n^-\ce\lleq \ga\ce\lleq \lambda_n^+\ce$. Then, applying $D$ is a dynamical realization-like homomorphism for $\lleq$, we deduce $D(\lambda_n^-)(p)\leq D(\ga)(p)\leq D(\lambda_n^+)(p)$, as we wanted. This finishes the proof of the first statement.\\

We now show the second statement. Consider $\ga\in \bola_{\equis_{\gr}}(n+1)$. This means that we can write $\ga=\alpha_1\cdots\alpha_j$ for some $\alpha_1,\dots,\alpha_j\in \equis_{\gr}$ and where $j\leq n+1$. Notice that $\alpha_i\cdots\alpha_1\in \bola_{\equis_{\gr}}(n)$ for all $i< j$, so we deduce that $\lambda_n^-\ce\lleq\alpha_i\cdots\alpha_1\ce\lleq \lambda_n^+\ce$ for all $i< j$. Then, applying $D$, we obtain that the elements $$x_0= p,\; x_1=D(\alpha_1)(p),\;\dots,\;x_{j-1}=D(\alpha_{j-1}\cdots\alpha_1)(p)$$ they all belong to $[D(\lambda_n^-)(p),D(\lambda_n^+)(p)]$. Hence, applying the hypothesis, we obtain recursively that 
\begin{equation*}
\begin{split}
\tilde{D}(\alpha_k\alpha_{k-1}\cdots\alpha_1)(p)=\tilde{D}(\alpha_k)\left(\tilde{D}(\alpha_{k-1}\cdots\alpha_1)(p)\right)=\tilde{D}(\alpha_k)(x_{k-1})=\\
= D(\alpha_k)(x_{k-1})=D(\alpha_k)\left(D(\alpha_{k-1}\cdots\alpha_1)(p)\right)=D(\alpha_k\alpha_{k-1}\cdots\alpha_1)(p)=x_k
\end{split}
\end{equation*}
for all $k=1,\dots, j$.
\end{proof}
\end{lema}

We can extend left-preorders induced from a dynamical realization. We state the needed property as a corollary of a more general that does not need to consider dynamical realizations.

\begin{prop}\label{extension orden}
Let $\gr$ be a group and $\ce,\ce_*$ be proper subgroups on $\gr$ such that $\ce\subseteq \ce_*$. Let $\lleq$ be a left-preorder on $\gr$ relative to $\ce_*$. Assume that $\ce$ is left-relatively convex on $\ce_*$. Then, there exists $\lleq'$ a left-preorder on $\gr$ relative to $\ce$ extending $\lleq$, meaning that $\ga_1\ce_*\lle \ga_2\ce_*$ implies $\ga_1\ce\lle' \ga_2\ce$.
\begin{proof}
Since $\ce$ a left-relatively convex subgroup on $\ce_*$, we obtain $\lleq_0$ a left-preorder on $\ce_*$ relative to $\ce$. Consider $P_0$ the positive cone of $\lleq_0$ and $P$ the positive cone of $\lleq$. Define $P'=P\cup P_0$. We now prove that $P'$ is the positive cone of a left-preorder $\lleq'$ on $\gr$ relative to $\ce$ by applying \Cref{orden relativo como semigrupo}.

Firstly, $P'$ is a semigroup. Indeed, let $\ga_1,\ga_2\in P'$. We have four cases. For the case $\ga_1,\ga_2\in P_0$, since $P_0$ is a subsemigroup as it is a positive cone, we have $\ga_1\ga_2\in P_0$. Then, in this case we have $\ga_1\ga_2\in P$. For the case $\ga_1,\ga_2\in P$, using that $P$ is a subsemigroup as it is a positive cone, we have $\ga_1\ga_2\in P$. Then, we obtain $\ga_1\ga_2\in P$. For the case $\ga_1\in P_0$ and $\ga_2\in P$, since $\ce P\subseteq P$ (recall \Cref{positivo por neutro es positivo}) and since $P_0\subseteq \ce$, we deduce that $\ga_1\ga_2$ belongs to $P\subseteq P'$. The case $\ga_1\in P$ and $\ga_2\in P_0$ is analogous to the previous case.\\

Secondly, we have $$\gr=P\sqcup P\inv \sqcup \ce_*=P\sqcup P\inv \sqcup P_0\sqcup P_0\inv \sqcup \ce=(P\sqcup P_0) \sqcup (P\sqcup P_0)\inv \sqcup \ce.$$ With this, we conclude that $P'$ is the positive cone of a left-preorder $\lleq'$ on $\gr$ relative to $\ce$.\\

Finally, we prove that $\lleq'$ satisfies the extension condition, that is, $\ga_1\ce_*\lle \ga_2\ce_*$ implies $\ga_1\ce\lle' \ga_2\ce$. This is direct from the fact that $P\subseteq P'$ by definition and from the fact that $P$ is the positive cone of $\lleq$ and $P'$ is the positive cone of $\lleq'$.
\end{proof}
\end{prop}

\begin{coro}\label{extension desde dinamica}
Let $\gr$ be a group and $\ce$ a left-relatively convex subgroup on $\gr$. Let $\df{D}{\gr}{\homeomr}$ be a homomorphism and let $p\in \real$ be such that $\ce\subseteq \stab_D(p)$ and $\stab_D(p)\neq \gr$. Consider $\lleq_{p}$ the relation on $\gr/\ce$ defined by $\ce\lleq_{p} \ga \ce \Leftrightarrow p\leq D(\ga)(p)$. Then $\lleq_{p}$ can be extended to a \defrelativeorder{\gr}{\ce}, say $\lleq$, in the sense that $\ga_1\ce\lle_{p} \ga_2\ce$ implies that $\ga_1\ce\lle \ga_2\ce$ for all $\ga_1,\ga_2\in\gr$.
\begin{proof}
When $\ce= \stab_D(p)$, the result follows directly by taking $\lle_{p}=\lleq$, so assume that $\ce\subsetneq \stab_D(p)$. Since $\ce$ is a left-relatively convex subgroup on $\gr$, there is some left-preorder on $\gr$ relative to $\ce$. Then, when restricting this left-preorder to $\stab_D(p)$, which is possible since $\ce\subseteq \stab_D(p)$, we obtain $\lleq^*$ a left-preorder on $\stab_D(p)$ relative to $\ce$. In other words, we have shown that $\ce$ is relatively convex on $\stab_D(p)$. Notice also that $\lleq_{p}$ is a left-preorder on $\gr$ relative to $\stab_D(p)$ by applying \Cref{lema realizacion dinamica orden}. Hence, the result follows by applying \Cref{extension orden}.
\end{proof}
\end{coro}

\section{The space of left-preorders}

We start by introducing the space of left-preorders. The definition of the space of left-preorders, as it is done in done in \cite{morris-antolinrivas}, follows the ideas of the space of left-orders considered by Sikora in \cite{topology-sikora}.

\begin{nota}
Let $\gr$ be a group and $\ce$ a left-relatively convex subgroup. By \Cref{orden relativo como semigrupo}, given $\lleq$ a left-preorder on $\gr$ relative to $\ce $ and denoting by $P$ its positive cone, we can define a $\phi_\lleq\in \{\pm1,0\}^\gr$ by:
\begin{equation*}
\begin{array}{r@{\hspace{0pt}}c@{\hspace{0pt}}
c@{\hspace{4pt}}l}
\phi_\lleq\colon &\gr &\longrightarrow& \{\pm1,0\} \\

&\ga \ \ &\mapsto&\begin{cases}1 \quad \textrm{if} \quad \ga \in P \\
 -1 \quad \textrm{if}\quad \ga \in P\inv \\
 0 \quad \textrm{if}\quad \ga \in \ce \end{cases}.
\end{array}
\end{equation*}
In particular, we can identify the set $\pord(\gr)$ of left-preorders on $\gr$ with a subset of $\{\pm1,0\}^\gr$. The \textit{space of left-preorders} on $\gr$ is the set $\pord(\gr)$ endowed with the subspace topology of $\{\pm 1,0\}^\gr$ with the product topology. The \textit{space of left-preorders on} $\gr$ \textit{relative to} $\ce$ is the set $\pord_\ce (\gr)$ endowed with the subspace topology of $\pord(\gr)$.
\end{nota}

\begin{obs}
A basis of open neighborhoods of $\phi\in\pord(\gr)$ is given by the sets $$U_{g_1,\dots,g_n}(\phi)=\{\psi\in\pord(\gr) \tq \psi(g_i)=\phi(g_i)\}$$ for all $g_1,\dots, g_n\in \gr$ and all $n\in \natp$. The sets $$\overline{U}_{g_1,\dots,g_n}(\phi)=\{\psi\in\pord(\gr) \tq \psi(g_i)\neq\phi(g_i)\}$$ are also open. Notice that these facts are also true when changing $\pord(\gr)$ by $\pord_\ce (\gr)$.
\end{obs}

We now study some different topological properties of these spaces. In some properties, we will sometimes consider also the space $\pord (\gr)\cup\{\lleq_{T}\}$ for $\lleq_{T}$ the trivial preorder. We do it only because it will be helpful to prove \Cref{corolario espacio de preordenes es incontable producto libre}. Essentially, it is sometimes convenient to consider $\pord (\gr)\cup\{\lleq_{T}\}$ instead of $\pord (\gr)$ to guarantee compactness without any extra condition in $\gr$.\\

The next two propositions follow directly since the properties considered are satisfied by the space $\{\pm1,0\}^\gr$ and they are hereditary for subspaces.

\begin{prop}\label{espacio de ordenes es totalmente disconexo y hausdorff}
Let $\gr$ be a group. The spaces $\pord(\gr)$, $\pord (\gr)\cup\{\lleq_{T}\}$ and $\pord_\ce (\gr)$ are totally disconnected and Hausdorff for $\lleq_{T}$ the trivial preorder and for $\ce$ any left-relatively convex subgroup.
\end{prop}

\begin{prop}\label{espacio de ordenes en contable es espacio metrico}
Let $\gr$ be a group. If $\gr$ is countable, then $\pord(\gr)$, $\pord (\gr)\cup\{\lleq_{T}\}$ and $\pord_\ce (\gr)$ are metrizable for $\lleq_{T}$ the trivial preorder and for $\ce$ any left-relatively convex subgroup.
\end{prop}

\begin{nota}
Let $\gr$ be a group. We define the following subsets of $\{\pm1,0\}^\gr$:
\begin{equation*}
V_{\ga }^{r}=\{\psi\in \{\pm1,0\}^\gr \tq \psi(\ga )=r\} \quad \textrm{for each} \; \ga \in \gr\; \textrm{and each }\; r\in\{\pm1,0\},
\end{equation*}
\begin{equation*}
\overline{V}_{\ga }^{r}=\{\psi\in \{\pm1,0\}^\gr\tq \psi(\ga )\neq r\}=\bigcup_{k\in\{\pm 1, 0\}\settminus\{r\}}V_{\ga }^{k} \quad \textrm{for each} \; \ga \in \gr\; \textrm{and each }\; r\in\{\pm1,0\}.
\end{equation*}
Notice that they are clopen (simultaneously open and closed) sets on $\{\pm1,0\}^\gr$.
\end{nota}

We prove now that $\pord_\ce (\gr)$ is compact by showing that it is closed on the compact space $\{\pm1,0\}^\gr$. The fact that $\{\pm1,0\}^\gr$ is compact is a consequence of Tychonoff's theorem.

\begin{prop}\label{espacio de ordenes relativos es cerrado}
Let $\gr$ be a group and $\ce$ a left-relatively convex subgroup of $\gr$. Then, the set $\pord_\ce (\gr)$ is closed on $\{\pm1,0\}^\gr$. In particular, $\pord_\ce (\gr)$ is a compact space.
\begin{proof}
We prove that $\pord_\ce (\gr)$ is closed by describing it as a complement of an open set made of intersections and finite unions of open sets. Consider the open sets $V_{\ga}^{r}$ and $\overline{V}_{\ga}^{r}$ of $\{\pm1,0\}^\gr$ as above. From \Cref{orden relativo como semigrupo} we deduce that having an element of $\pord_\ce (\gr)$ is equivalent to give a decomposition of the form $\gr=P\sqcup P\inv \sqcup \ce $ where $P$ is a subsemigroup of $\gr$. This is equivalent to give an element $\phi$ of $\{\pm1,0\}^\gr$ satisfying the following conditions:
\begin{enumerate}
\item $\left(\phi\inv(\{1\})\right)\inv=\phi\inv(\{-1\})$.
\item $\phi\inv(\{0\})=\ce $.
\item $\phi\inv(\{1\})$ is a subsemigroup of $\gr$.
\end{enumerate}
This give us the description
\begin{equation*}
\begin{split}
\{\pm1,0\}^\gr \settminus \pord_\ce (\gr)=&\left((\bigcup_{\ga \in \gr}V_{\ga \inv}^{1}\cap \overline{V}_{\ga }^{-1}))\cup (\bigcup_{\ga \in \gr}\overline{V}_{\ga \inv}^{1}\cap V_{\ga }^{-1}))\right)\cup \left((\bigcup_{\ga \in \ce }\overline{V}_{\ga }^{0}))\cup (\bigcup_{\ga \in \ce \settminus \gr} V_{\ga }^{0}))\right)\cup\\ & \cup\left(\bigcup_{\ga ,\gb \in \gr}V_{\ga }^{1}\cap V_{\gb }^{1}\cap \overline{V}_{\ga \gb }^{1}\right)
\end{split}
\end{equation*}
where each of the four sets between parentheses is an open set that corresponds to the respective negation of each condition below. Therefore, we have shown that $\pord_\ce (\gr)$ is a closed subset.
\end{proof}
\end{prop}

For the compactness of $\pord (\gr)$ we need to assume that $\gr$ is finitely generated. However, $\pord (\gr)\cup\{\lleq_{T}\}$ will be compact in all cases, for $\lleq_{T}$ be the trivial preorder.

\begin{prop}\label{espacio de ordenes con el trivial es cerrado}
Let $\gr$ be a group. Let $\lleq_{T}$ be the trivial preorder. Then, the set $\pord (\gr)\cup\{\lleq_{T}\}$ is closed on $\{\pm1,0\}^\gr$. In particular, $\pord (\gr)\cup\{\lleq_{T}\}$ is a compact space.
\begin{proof}
The proof is similar to the one done in \Cref{espacio de ordenes relativos es cerrado}. Consider the open sets $V_{\ga}^{r}$ and $\overline{V}_{\ga}^{r}$ of $\{\pm1,0\}^\gr$ as above. By \Cref{orden relativo como semigrupo}, having an element of $\pord (\gr)\cup\{\lleq_{T}\}$ is equivalent to give a decomposition of the form $\gr=P\sqcup P\inv \sqcup \ce $ where $P$ is a subsemigroup of $\gr$ satisfying that $\ce $ is a subgroup (for $\ce=\gr$ it give us $\lleq_{T}$). This is equivalent to give an element $\phi\in \{\pm1,0\}^\gr$ satisfying the following:
\begin{enumerate}
\item $\left(\phi\inv(\{1\})\right)\inv=\phi\inv(\{-1\})$.
\item $\phi\inv(\{0\})$ is a subgroup.
\item $\phi\inv(\{1\})$ is a subsemigroup of $\gr$.
\end{enumerate}
This give us the description
\begin{equation*}
\begin{split}
\{\pm1,0\}^\gr \settminus \left(\pord (\gr)\cup\{\lleq_{T}\}\right)=&\left((\bigcup_{\ga \in \gr}V_{\ga \inv}^{1}\cap \overline{V}_{\ga }^{-1}))\cup (\bigcup_{\ga \in \gr}\overline{V}_{\ga \inv}^{1}\cap V_{\ga }^{-1}))\right)\cup\left(\bigcup_{\ga ,\gb \in \gr}V_{\ga }^{0}\cap V_{\gb }^{0}\cap \overline{V}_{\ga \gb\inv }^{0}\right) \cup\\ & \cup\left(\bigcup_{\ga ,\gb \in \gr}V_{\ga }^{1}\cap V_{\gb }^{1}\cap \overline{V}_{\ga \gb }^{1}\right)
\end{split}
\end{equation*}
Therefore, $\pord(\gr)\cup\{\lleq_{T}\}$ is closed on $\{\pm1,0\}^\gr$.
\end{proof}
\end{prop}

\begin{prop}\label{espacio de ordenes es cerrado}
Let $\gr$ be a group. Assume that $\gr$ is finitely generated. Then, the set $\pord (\gr)$ is closed on $\{\pm1,0\}^\gr$. In particular, $\pord (\gr)$ is a compact space.
\begin{proof}
Since $\gr$ is finitely generated, we can take $\gr=\gen{\ga_1,\dots,\ga_n}$. Consider the open sets $V_{\ga}^{r}$ as before. Let $\lleq_{T}$ be the trivial preorder. By \Cref{orden relativo como semigrupo}, having an element of $\pord (\gr)\cup\{\lleq_{T}\}$ is equivalent to give a decomposition of the form $\gr=P\sqcup P\inv \sqcup \ce $ where $P$ is a subsemigroup of $\gr$ satisfying that $\ce $ is a subgroup (for $\ce=\gr$ it give us $\lleq_{T}$). In particular, since $\gr=\gen{\ga_1,\dots,\ga_n}$, we have $$\left(\pord (\gr)\cup\{\lleq_{T}\}\right)\cap\left(\bigcap_{i=1, \dots, n}V_{\ga_i}^{0}  \right)=\{\lleq_{T}\}.$$ Then, since $V_{\ga_i}^{0}$ is open on $\{\pm1,0\}^\gr$, we deduce that $\{\lleq_{T}\}$ is open on $\pord (\gr)\cup\{\lleq_{T}\}$. In particular, $\pord (\gr)$ is closed on $\pord (\gr)\cup\{\lleq_{T}\}$. Since $\pord (\gr)\cup\{\lleq_{T}\}$ is closed on $\{\pm1,0\}^\gr$ by \Cref{espacio de ordenes con el trivial es cerrado}, we conclude that $\pord (\gr)$ is closed on $\{\pm1,0\}^\gr$.
\end{proof}
\end{prop}

We need the condition of $\gr$ being finitely generated to guarantee that $\pord (\gr)$ is compact. For example, $\pord \left(\bigoplus_{i\in \nat}\ent \right)$ is not compact, as it is shown in Example 2.4 of \cite{morris-antolinrivas}.\\

In some situations these spaces will become a Cantor set. In this sense, it is very useful the following topological characterization since the spaces we consider are compact, Hausdorff, totally disconnected, and non-empty metric spaces (depending on the case, we may need to require that $\gr$ is finitely generated or countable, as discussed in the previous results). For a proof, see for instance Corollary 2-98 from \cite{topology-hockingyoung}.

\begin{prop}\label{cuando es cantor}
Every compact, Hausdorff, totally disconnected, non-empty metric space having no isolated points is homeomorphic to the Cantor set.
\end{prop}

We now establish a relation between being isolated in the space of left-preorders and having a positive cone that is finitely generated as a semigroup. This behaviour is well known for left-orders, due to Linnell.

\begin{nota}
Given $\gr$ a group and $S\subset \gr$, $\semgen{S}$ denotes the subsemigroup generated by $S$.
\end{nota}

\begin{prop}\label{fg como semigrupo implica asilado}
Let $\gr$ be a group and $\ce$ a left-relatively convex subgroup of $\gr$. Let $\lleq$ be a \defrelativeorder{\gr}{\ce} having positive cone $P$. If $P$ is finitely generated as a subsemigroup, then $\lleq$ is isolated in $\pord_\ce(\gr)$. Moreover, if we also assume that $\ce$ is a finitely generated subgroup, then $\lleq$ is also isolated in $\pord(\gr)$.
\begin{proof}
Consider $P=\semgen{\ga_1,\dots,\ga_n}$ the positive cone of $\lleq$, so that $\gr=P\sqcup P\inv \sqcup \ce$. In order to prove that $\lleq$ is isolated in $\pord_\ce(\gr)$ we need to prove that $\{\lleq\}$ is an open set. For this purpose, we are going to prove that $$\{\lleq'\in \pord_\ce(\gr)\tq \phi_{\lleq'}(\ga_1)=1, \dots,\phi_{\lleq'}(\ga_n)=1\}=\{\lleq\},$$ so the fact that $\{\lleq\}$ is an open set will follow directly. Indeed, take $\lleq'$ an element belonging to the set in the left term of the equation. Then, $\lleq'$ is a left-preorder on $\gr$ relative to $\ce$ such that $\ce\lleq'\ga_i \ce$ for $i=1,\dots,n$. Denote by $P'$ its positive cone. Using the fact that $P=\semgen{\ga_1,\dots,\ga_n}$, we obtain that $$P\subseteq P'\mbox{ and }P\inv\subseteq (P')\inv.$$ Notice also that $\gr=P'\sqcup (P')\inv \sqcup \ce$ by \Cref{orden relativo como semigrupo}. Hence, $P=P'$. We deduce that $\lleq'=\lleq$, as we wanted to prove. \\

Assume now that $\ce$ is a finitely generated subgroup. Considering $\gx_1,\dots,\gx_k$ a finite set of generators, following similar arguments as before we can conclude that $$\{\lleq'\in \pord(\gr)\tq \phi_{\lleq'}(\ga_1)=1, \dots,\phi_{\lleq'}(\ga_n)=1, \phi_{\lleq'}(\gx_1)=0, \dots,\phi_{\lleq'}(\gx_k)=0 \}=\{\lleq\}.$$ Therefore, $\lleq$ is isolated in $\pord(\gr)$.
\end{proof}
\end{prop}

We can use compactness of $\{\pm 1,0\}^G$ to obtain a way to left-preorder a whole group by just left-preordering all finitely generated subgroups. We do it in the next proposition. This will help us to prove \Cref{teorema aislados} in such a way we will be able to restrict to the finitely generated case. Notice that, in next defintion and next proposition, $C$ is not necessarily going to be left-relatively convex on $G$.
\begin{defi}
Let $G$ be a group and $C$ be a proper subgroup of $G$. Take $x_1,\dots,x_n\in G\settminus C$ where $n\in\natp$. We say that $\{x_1,\dots,x_n\}$ has property $(E_C)$ if and only if for every finite family $f_1,\dots,f_k$ of elements of $G\settminus C$ there exists $\eta_1\dots,\eta_k\in \{\pm 1\}$ such that $C\cap \semgen{x_1,\dots,x_n,f_1^{\eta_1},\dots,f_k^{\eta_k}}=\emptyset$, where $\semgen{S}$ is the subsemigroup generated by $S$. Such $\eta_i$ exponents are called compatible.
\end{defi}

\begin{prop}\label{propiedad E}
Let $G$ be a group and $C$ be a proper subgroup of $G$. Take $x_1,\dots,x_n\in G\settminus C$ where $n\in\natp$. Then, $G$ admits a left-preorder $\lleq$ relative to $C$ such that $C\lle x_i C$ for all $i=1\dots n$ if and only if $\{x_1,\dots,x_n\}$ has property $(E_C)$.
\begin{proof}
The necessity of the equivalence is clear. We now prove the other implication using topological properties. Assume $\{x_1,\dots,x_n\}$ has property $(E_C)$. Consider the simultaneusly open and closed sets $V_{g}^{r}$ and $\overline{V}_{g}^{r}$ of $\{\pm1,0\}^G$ as above. 
We also define the following closed sets of $\{\pm 1,0\}^G$:
\begin{equation*}
\begin{split}
A_0=&\bigcap_{x\in G}\left((V_{x}^{1}\cap V_{x\inv}^{-1})\cup (V_{x}^{-1}\cap V_{x\inv}^{1})\cup (V_{x}^{0}\cap V_{x\inv}^{0})\right)\\
A_1=&\{\pm 1,0\}^G\settminus\bigcup_{x,y\in G}\left(V_{x}^{1}\cap V_{y}^{1}\cap \overline{V}_{xy}^{1}\right)\\
A_2=&\bigcap_{g\in C}V_{g}^{0}\\
A_3=&\bigcup_{i=1,\dots, n}V_{x_i}^{1}
\end{split}
\end{equation*}
Notice that the following equivalences are direct consequence of the definitions above:
\begin{equation*}
\begin{split}
\chi\in A_0 \sii& \left(\chi\inv(\{1\})\right)\inv=\chi\inv(\{-1\})\\
\chi\in A_1 \sii& \chi\inv(\{1\})\; \textrm{ is a subsemigroup of } G\\
\chi\in A_2 \sii& C\subseteq \chi\inv(\{0\})\\
\chi\in A_3 \sii& x_1,\dots,x_n\in\chi\inv(\{1\})
\end{split}
\end{equation*}
We define for each $f_1,\dots,f_l\in G\settminus C$ and each $l\in\natp$ the set
\begin{equation*}
\tilde{A}_{\{f_1,\dots,f_l\}}=(V_{f_1}^{1}\cup V_{f_1\inv}^{1})\cap\dots \cap (V_{f_n}^{1}\cup V_{f_n\inv}^{1}) \quad \textrm{for each} \; f_1,\dots,f_l\in G\settminus C\; \textrm{ and }\; l\in\natp.
\end{equation*}
In this way, we have
$$\chi\in \tilde{A}_{\{f_1,\dots,f_l\}} \sii f_1^{\epsilon_1},\dots,f_l^{\epsilon_l}\in\chi\inv(\{1\}) \; \textrm{ for some }\; \epsilon_1,\dots,\epsilon_l\in \{\pm 1\}.$$
We also define for each $f_1,\dots,f_l\in G\settminus C$ and each $l\in\natp$ the set
\begin{equation}\label{propiedad E eq1}
B_{\{f_1,\dots,f_l\}}=A_0\cap A_1\cap A_2\cap A_3\cap \tilde{A}_{\{f_1,\dots,f_l\}}\quad \textrm{for each} \; f_1,\dots,f_l\in G\settminus C\; \textrm{ and }\; l\in\natp.
\end{equation}

We are now going to prove that
\begin{equation}\label{propiedad E eq2}
B_{\{f_1,\dots,f_l\}}\neq \emptyset  \;\mbox{ for all }f_1,\dots,f_l\in G\settminus C \mbox{ where }l\in\natp.
\end{equation}
Indeed, consider $f_1,\dots,f_l\in G\settminus C$ where $l\in\natp$. Since $\{x_1,\dots,x_n\}$ has property $(E_C)$ by hypothesis, we know that there are $\eta_1\dots,\eta_l\in \{\pm 1\}$ such that $C\cap S=\emptyset$ where $S=\semgen{x_1,\dots,x_n,f_1^{\eta_1},\dots,f_l^{\eta_l}}$. Notice that $S\cap S\inv=\emptyset$ as, otherwise, it would imply $1\in S$, which is a contradiction since $C\cap S=\emptyset$. Therefore, we can define $\chi\in\{\pm 1, 0\}^G$ by $$\chi\inv(\{1\})=S\mbox{ , }\chi\inv(\{-1\})=S\inv \mbox{ and } \chi\inv(\{0\})=G\settminus (S\cup S\inv).$$ We are going to prove that $$\chi\in B_{\{f_1,\dots,f_l\}}.$$ By \Cref{propiedad E eq1}, we only need to prove that $$\chi\in A_0\mbox{ , }\chi\in A_1\mbox{ , }\chi\in A_2\mbox{ , }\chi\in A_3\mbox{ and }\chi\in \tilde{A}_{\{f_1,\dots,f_l\}}$$ It is clear that $\chi\in A_0$. Also $\chi\in A_1$ since $S$ is a subsemigroup of $G$. We know that $C\cap S=\emptyset$, but taking inverses we can also deduce that $C\cap S\inv=\emptyset$, so then $C\subseteq G\settminus (S\cup S\inv)=\chi\inv(\{0\})$, and hence $\chi\in A_2$. By definition of $S$ we have $x_1,\dots,x_n\in S=\chi\inv(\{1\})$, so $\chi\in A_3$. Also by definition of $S$ it is clear that $f_1^{\eta_1},\dots,f_l^{\eta_l}\in S=\chi\inv(\{1\})$, so $\chi\in \tilde{A}_{\{f_1,\dots,f_l\}}$. Hence, we have shown that $$\chi\in B_{\{f_1,\dots,f_l\}},$$ and this proves \Cref{propiedad E eq2}, as we wanted. \\

We can consider the family $\Omega$ of non-empty closed subsets of $\{\pm 1, 0\}^G$ given by the sets $B_{\{f_1,\dots,f_l\}}$ for $f_1,\dots,f_l\in G\settminus C$ and $l\in\natp$. We claim that 
\begin{equation}\label{propiedad E eq3}
\bigcap_{B\in \Omega}B\neq \emptyset.
\end{equation}
Since all the elements in $\Omega$ are closed sets of the compact space $\{\pm 1, 0\}^G$, \Cref{propiedad E eq3} will follow from the fact that $\Omega$ has the finite intersection property. We now prove that $\Omega$ has the finite intersection property. Consider $B_{\{f_{1,1},\dots,f_{1,l_1}\}}, \dots, B_{\{f_{t,1},\dots,f_{t,l_t}\}}\in \Omega$. By following the definitions, we have $$B_{\{f_{1,1},\dots,f_{1,l_1}, \dots, \dots,f_{t,1},\dots,f_{t,l_t}\}}\subseteq B_{\{f_{1,1},\dots,f_{1,l_1}\}}\cap \dots\cap B_{\{f_{t,1},\dots,f_{t,l_t}\}},$$ but $B_{\{f_{1,1},\dots,f_{1,l_1}, \dots \dots,f_{t,1},\dots,f_{t,l_t}\}}$ is non-empty by \Cref{propiedad E eq2}, so the intersection above is non-empty. Then, we have shown that $\Omega$ has the finite intersection property, so \Cref{propiedad E eq3} holds.\\

We are now going to prove that there exists a $\chi\in \{\pm 1, 0\}^G$ that will give us a left-preorder with the needed properties. By \Cref{propiedad E eq3}, we can consider $$\chi\in\bigcap_{B\in \Omega}B.$$ By definition of $\Omega$, we have $\chi\in B_{\{f_1,\dots,f_l\}}$ for all $f_1,\dots,f_l\in G\settminus C$ and all $l\in\natp$. By \Cref{propiedad E eq1}, we have that $\chi\in A_i$ for $i=0,1,2,3$ and that $\chi\in \tilde{A}_{\{f_1,\dots,f_l\}}$ for all $f_1,\dots,f_l\in G\settminus C$ and all $l\in\natp$. We claim that 
\begin{equation}\label{propiedad E eq4}
C= \chi\inv(\{0\}).
\end{equation}
We are going to prove it. Since $\chi\in A_2$, we deduce that $C\subseteq \chi\inv(\{0\})$. Then, we only need to prove $$\chi\inv(\{0\})\subseteq C.$$ Indeed, given $f\in G\settminus C$, we know that $\chi\in \tilde{A}_{\{f\}}$, so $f^{\epsilon}\in\chi\inv(\{1\})$ for some $\epsilon\in \{\pm 1\}$. Notice that $\left(\chi\inv(\{1\})\right)\inv=\chi\inv(\{-1\})$ since $\chi\in A_0$, so we deduce $f\in\chi\inv(\{\pm 1\})$. Then, $f\notin\chi\inv(\{0\})$. Hence, we have shown \Cref{propiedad E eq4}.\\

We are now ready use $\chi$ to obtain the left-preorder with the needed properties. We can decompose $G$ in a disjoint union as $$G=\chi\inv(\{1\})\sqcup\chi\inv(\{-1\})\sqcup\chi\inv(\{0\}).$$ Define $P=\chi\inv(\{1\})$. Notice that $P$ is a subsemigroup of $G$ since $\chi\in A_1$. Notice also that $P\inv=\chi\inv(\{-1\})$ as a consequence of $\chi\in A_0$. Then, we can write the the decomposition bellow as $$G=P\sqcup P\inv \sqcup \chi\inv(\{0\})$$ where $P$ is a subsemigroup of $G$. Therefore, using \Cref{propiedad E eq4}, the decomposition bellow is in fact $$G=P\sqcup P\inv \sqcup C$$ where $P$ is a subsemigroup of $G$. Hence, by \Cref{orden relativo como semigrupo}, there is a left-preorder $\lleq$ on $G$ relative to $C$ with positive cone $P$. Finally, notice that $x_1,\dots x_n\in P$ because $\chi\in A_3$. We conclude that $G$ admits a left-preorder $\lleq$ relative to $C$ such that $C\lle x_i C$ for all $i=1\dots n$.
\end{proof}
\end{prop}

We conclude the section by giving as an example a brief description of some properties of the space of left-preorders for $\ent^m$. This description is done in \cite{valuated-decauprond} for abelian groups with all the details. In this paper, the notion of Patch topology on the set of left-preorders corresponds to the topology we consider on the space of left-preorders. 

\begin{ej}
Consider the group $\ent^m$ for $m\in\natp$. Consider $C$ a (normal) subgroup $C$ on $\ent^m$, so it is free abelian of some rank $l\leq m$. If $C$ is left-relatively convex, then $G/C$ is left-orderable and abelian, so it is a free abelian group of rank $m-l>0$ and co-rank $l<m$.\\

The space $\pord_C(\ent^m)$ can be identified with the space of left-orderings of the group $\ent^m/C$, which is isomorphic to $\ent^{m-l}$. By \cite{topology-sikora}, where we can find a description of the space of left-orders of $\ent^r$, we deduce that $\pord_C(\ent^m)$ is a two point discrete space if $C$ has abelian co-rank $1$ and a Cantor set if $C$ has abelian co-rank greater than $1$.\\

We can now describe $\pord(\ent^m)$. For $m=1$, $\pord(\ent)$ is a two point discrete space. Assume $m\geq 2$. We will show that $\pord(\ent^m)$ consists of a numerable discrete space and a Cantor set. Consider $C$ a left-relatively convex subgroup of abelian co-rank $1$, so we can consider $c_1,\dots,c_{m-1}$ a set of generators of $C$ and $xC$ a generator of the cyclic group $\ent^m/C$. Since $\ent^m/C$ is cyclic, there are only $2$ left-preorders on $\ent^m$ relative to $C$, say $\lleq_C^+$ and $\lleq_C^-$. Then, we have $$\{\lleq'\in \pord(\ent)\tq \phi_{\lleq}(c_1)=0,\dots,\phi_{\lleq}(c_{m-1})=0,\phi_{\lleq}(x)=1 \}=\{\lleq_C^+\} \;\mbox{ and }$$ $$\{\lleq'\in \pord(\ent)\tq \phi_{\lleq}(c_1)=0,\dots,\phi_{\lleq}(c_{m-1})=0,\phi_{\lleq}(x)=-1 \}=\{\lleq_C^-\}.$$ Hence, both left-preorders are isolated.\\

On the other hand, assume $C$ is a left-relatively convex subgroup of abelian co-rank greater than $1$. Take $\lleq$ a left-preorder on $\ent^m/C$. Then $\lleq$ is not isolated on $\pord_C(\ent^m)$ as it is a Cantor set. In particular, $\lleq$ is not isolated on $\pord(\ent^m)$. Hence, the set of left-preorders relative to a subgroup of abelian co-rank greater than $1$ has no isolated elements. In fact, it is a Cantor set, as it is a closed subset of $\pord(\ent^m)$.\\

We conclude that $\pord(\ent^m)$ consists of a numerable discrete space and a Cantor set. The numerable discrete space is the set of left-preorders relative to a subgroup of abelian co-rank $1$ and the Cantor set is the set of left-preorders relative to a subgroup of abelian co-rank greater than $1$.
\end{ej}

\section{The space of left-preorders on free products}

We now introduce the definition of the Kurosh rank of a free product. We use the definition given in \cite{burns-kurosh}. We previously need to state Kurosh's theorem about the structure of subgroups of a free product (for a reference see Theorem 14 of \cite{serre}).

\begin{teo}[Kurosh's theorem]
Let $\gr,\grb$ be groups. Let $\ce$ be a subgroup of $\gr*\grb$. Then, 
\begin{equation}\label{eq teorema kurosh}
\ce=\mathbb{F}*\left(\ast \ce\cap \gx\gr \gx\inv\right)*\left(\ast \ce\cap \gy\grb \gy\inv\right)
\end{equation} where $\mathbb{F}$ is a free subgroup of $\gr*\grb$ and where $\gx$ ranges over a set of representatives of the double coset $\ce\backslash(\gr*\grb)/\gr$ and $\gy$ ranges over a set of representatives of the double coset $\ce\backslash(\gr*\grb)/\grb$.
\end{teo}

\begin{defi}
Let $\gr,\grb$ be groups. Let $\ce$ be a subgroup of $\gr*\grb$. When considering a decomposition as in \Cref{eq teorema kurosh}, we define the \textit{Kurosh rank} of $\ce$ with respect to the free product $\gr*\grb$ as the free rank of $\mathbb{F}$ plus the number of non trivial factors $\ce\cap \gx\gr \gx\inv$ and $\ce\cap \gy\grb \gy\inv$.
\end{defi}

Notice that the Kurosh rank is possibly non-finite. This notion is well defined and does not depend on the chosen representatives (see \cite{burns-kurosh}). Every finitely generated subgroup has finite Kurosh rank.\\

\begin{obs}\label{obs bs rango kurosh infinito}
Let $\gr,\grb$ be groups. Let $\ce$ be a subgroup of $\gr*\grb$. If $\ce$ has infinite Kurosh rank with respect to $\gr*\grb$, then for each finite subset of $\ce$ there exists a proper free factor $\ce_1$ of $\ce$ containing that subset.
\end{obs}

We can now focus on the main results about the spaces $\pord(\gr*\grb)$ and $\pord_\ce(\gr*\grb)$. We state the theorem about $\pord_\ce(\gr*\grb)$. This result will be used in order to prove the result about $\pord(\gr*\grb)$.

\begin{teo}\label{teorema aislados}
Let $\gr$ and $\grb$ be two non-trivial groups. Let $\ce$ be a left-relatively convex subgroup on $\gr*\grb$. If $\ce$ has finite Kurosh rank with respect to $\gr*\grb$, then $\pord_\ce(\gr*\grb)$ has no isolated elements.
\end{teo}

The theorem will follow from the following result, which allows us to reduce to the case where $\gr$ and $\grb$ are finitely generated.

\begin{teo}\label{teorema auxiliar aislados}
Let $\gr$ and $\grb$ be two non-trivial groups. Let $\lleq$ be a left-preorder on $\gr*\grb$ relative to a subgroup $\ce$ of finite Kurosh rank with respect to $\gr*\grb$. Let $\conjfiny$ be a finite subset of $\gr*\grb$ such that $\ce\lle \elefiny\ce$ for all $\elefiny\in\conjfiny$. Then, there exists finite sets $\equisg\subseteq \gr$ and $\equish\subseteq \grb$ and an element $\gamma_{*}\in \gen{\equisg,\equish}$ (they depend on the set $\conjfiny$) such that $\ce\lle \gamma_{*}\ce$  and satisfying the following property:\\

For all finitely generated subgroups $\gr_1\subseteq \gr$ and $\grb_1\subseteq \grb$ containing $\equisg$ and $\equish$ respectively, the group $\Gamma=\gr_1*\grb_1$ admits a left-preorder $\lleq^*$ relative to $\ce_\Gamma=\Gamma\cap \ce$ such that $\gamma_{*}\ce_\Gamma\lle^* \ce_\Gamma$ and $\ce_\Gamma\lle^* \elefiny\ce_\Gamma$ for all $\elefiny\in\conjfiny$.
\end{teo}

We now prove \Cref{teorema aislados} assuming \Cref{teorema auxiliar aislados}.

\begin{proof}[Proof of \Cref{teorema aislados}.]
Consider $\lleq$ a left-preorder on $\gr*\grb$ relative to $\ce$ with $\ce$ of finite Kurosh rank with respect to $\gr*\grb$. We need to show that $\lleq$ is not isolated on $\pord_\ce(\gr*\grb)$. In order to do that, we fix $\conjfiny$ a finite set of $\gr*\grb$, and we want to build another left-preorder on $\gr*\grb$ relative to $\ce$ that agrees on $\conjfiny$ with $\lleq$ but that is not equal to $\lleq$. Since all left-preorder on $\gr*\grb$ relative to $\ce$ agrees on $\ce$ with $\lleq$, we can assume that $\conjfiny$ has no elements in $\ce$. Then, by possibly changing some elements of $\conjfiny$ by its inverses, we can assume that all the elements $\elefiny$ in $\conjfiny$ satisfies that $\ce\lle \elefiny\ce$ or, in other words, we can assume that all the elements in $\conjfiny$ are $\lleq$-positive. Consider $\equisg\subseteq \gr$, $\equish\subseteq \grb$ and $\gamma_{*}\in \gen{\equisg,\equish}$ the sets and the element provided by \Cref{teorema auxiliar aislados} for the finite set $\conjfiny$. In order to build the needed new left-preorder we are going to use \Cref{propiedad E} with the set $\conjfiny\cup\{ \gamma_*\inv\}$, which is a finite set contained in $(\gr*\grb)\settminus \ce$. \\

To be able to use \Cref{propiedad E}, we are going to prove that $\conjfiny\cup\{ \gamma_*\inv\}$ satisfies property $(E_\ce)$. Take $\conjfinf=\{\elefinf_1,\dots,\elefinf_k\}$ a finite subset of $(\gr*\grb)\settminus \ce$. We need to prove that there exists $\eta_1,\dots,\eta_k\in \{\pm 1\}$ such that
\begin{equation}\label{teorema aislados eq1}
\semgen{\conjfiny\cup \{\gamma_*\inv\}\cup\{\elefinf_1^{\eta_1},\dots,\elefinf_k^{\eta_k} \}}\cap \ce=\emptyset.
\end{equation}
In order to do that, consider $\gr_1\subseteq \gr$ and $\grb_1\subseteq \grb$ two finitely generated subgroups such that $$\conjfiny\cup \{\gamma_*\inv\}\cup \equisg\cup \equish\cup \conjfinf\subset \Gamma=\gr_1*\grb_1.$$ By \Cref{teorema auxiliar aislados}, we know that there is a left-preorder $\lleq^*$ on $\Gamma$ relative to $\ce_\Gamma=\Gamma\cap \ce$ such that $\gamma_{*}\ce_\Gamma\lle^* \ce_\Gamma$ and $\ce_\Gamma\lle^* \elefiny\ce_\Gamma$ for all $\elefiny\in\conjfiny$. Since $\conjfinf\cap \ce$ is empty, for each $i=1,\dots,k$ we can take an $\eta_i\in \{\pm 1\}$ such that $\ce_\Gamma\lle^* \elefinf_i^{\eta_i}\ce_\Gamma$. Hence, we know that $\gamma_{*}\ce_\Gamma\lle^* \ce_\Gamma$, that $\ce_\Gamma\lle^* \elefiny\ce_\Gamma$ for all $\elefiny\in\conjfiny$ and that $\ce_\Gamma\lle^* \elefinf_i^{\eta_i}\ce_\Gamma$ for all $i=1,\dots,k$. Then, we deduce that the subsemigroup $\semgen{\conjfiny\cup \{\gamma_*\inv\}\cup\{\elefinf_1^{\eta_1},\dots,\elefinf_k^{\eta_k} \}}$ is contained in the positive cone of $\lleq^*$. By \Cref{orden relativo como semigrupo}, we deduce that $$\semgen{\conjfiny\cup \{\gamma_*\inv\}\cup\{\elefinf_1^{\eta_1},\dots,\elefinf_k^{\eta_k} \}}\cap (\Gamma\cap \ce)=\emptyset.$$ Since $\semgen{\conjfiny\cup \{\gamma_*\inv\}\cup\{\elefinf_1^{\eta_1},\dots,\elefinf_k^{\eta_k} \}}$ is contained in $\Gamma=\gr_1*\grb_1$ by definition of $\gr_1$ and $\grb_1$, we obtain \Cref{teorema aislados eq1}, as we wanted. Therefore, we have proven that $\conjfiny\cup\{ \gamma_*\inv\}$ satisfies property $(E_\ce)$.\\

As we have shown that $\conjfiny\cup\{ \gamma_*\inv\}$ satisfies property $(E_\ce)$, by \Cref{propiedad E}, there is a left-preorder $\hat{\lleq}$ on $\gr*\grb$ relative to $\ce$ such that $\ce\lle \gamma_*\inv \ce$ and $\ce\lle \elefiny\ce$ for all $\elefiny\in\conjfiny$. In particular, $\hat{\lleq}$ coincides with $\lleq$ on $\conjfiny$ but they are not equal since they not coincide on $\gamma_*$. Hence, we have built another left-preorder on $\gr*\grb$ relative to $\ce$ that agrees on $\conjfiny$ with $\lleq$ but is not equal to $\lleq$. We conclude that $\lleq$ is not isolated on $\pord_\ce(\gr*\grb)$, as we wanted.
\end{proof}

The proof of \Cref{teorema auxiliar aislados} uses dynamical properties following the ideas on \cite{freeproducts-rivas}.

\begin{proof}[Proof of \Cref{teorema auxiliar aislados}.]
Let $\lleq$ be a left-preorder on $\gr*\grb$ relative to a subgroup $\ce$ of finite Kurosh rank with respect to $\gr*\grb$. Let $\conjfiny$ be a finite subset of $\gr*\grb$ such that $\ce\lle \elefiny\ce$ for all $\elefiny\in\conjfiny$.\\

The fact that $\ce$ has finite Kurosh implies two facts that we are going to prove in next section. These two facts are the following:
\begin{itemize}
\item[\textbf{Fact I:}] There exists $\fino,\{\zeta_i\}_{i=1}^r,\{\alpha_j\}_{j=1}^s,\{\beta_t\}_{t=1}^l\subseteq \gr*\grb$ finite sets satisfying the following property: for every $\gr_0$ subgroup of $\gr$ and every $\grb_0$ subgroup of $\grb$ satisfying $\fino\subseteq \gen{\gr_0,\grb_0}$ there exists $\grfaca_j$ subgroup of $\gr_0$ for each $j=1,\dots,s$ and there exists $\grfacb_t$ subgroup of $\grb_0$ for each $t=1,\dots,l$ such that 
\begin{equation*}
\ce\cap\gen{\gr_0,\grb_0}=\gen{\{\zeta_i \tq i=1,\dots,r\}\cup \left(\bigcup_{j=1}^s\alpha_j\grfaca_j\alpha_j\inv \right)\cup \left(\bigcup_{l=1}^t\beta_l\grfacb_l\beta_l\inv \right)}.
\end{equation*}
\item[\textbf{Fact II:}] Let $(\lambda_l)_{l\in\nat}$ be a sequence of elements in $\gr*\grb$ such that $\lambda_l \ce\neq\lambda_m \ce$ for all $l\neq m$. Then, there exists $n\in\nat$ such that for all $k\geq n$ there are $\ga\in \equisg$ and $\gb\in \equish$ satisfying $\ga^{\varepsilon_\gr} \lambda_k \ce\neq\lambda_k \ce$, $\gb^{\varepsilon_\grb} \lambda_k \ce\neq\lambda_k \ce$ and $\ga^{\varepsilon_\gr} \lambda_k \ce\neq \gb^{\varepsilon_\grb}\lambda_k \ce$ for every $\varepsilon_\gr,\varepsilon_\grb\in\{\pm 1\}$.
\end{itemize}

Fix $\fino,\{\zeta_i\}_{i=1}^r,\{\alpha_j\}_{j=1}^s,\{\beta_t\}_{t=1}^l$ the elemnts given by Fact I. Consider $\equisg\subseteq \gr$ and $\equish\subseteq \grb$ two finite subsets, both having some non-trivial element, satisfying that $\equisg\inv=\equisg$ and $\equish\inv=\equish$ and satisfying that $$\conjfiny\cup\fino\subseteq \gen{\equisg}*\gen{\equish}.$$ For each $n\in\nat$, consider $\lambda_n^+$ a $\lleq$-maximal element in the ball $\bola_{\equis_{\gr}}(n)$ and consider $\lambda_n^-$ a $\lleq$-minimal element in the ball $\bola_{\equis_{\gr}}(n)$.\\

By \Cref{lema comparando lambdas}, we have $\lambda_n^+\ce\neq\lambda_k^+\ce$ for all $n,k\in\nat$. This implies that the sequence $(\lambda_n)_{n\in\nat}$ satisfies the hypothesis of Fact II. Using the finiteness of $\conjfiny$ and applying Fact II, we can consider a big enough $n\in\nat$ and we can take $\ga\in \equisg$ and $\gb\in \equish$ satisfying the following:
\begin{equation}\label{teorema auxiliar aislados eq1 x en bola}
\conjfiny\subseteq \bola_{\equisg\cup \equish}(n),
\end{equation}
\begin{equation}\label{teorema auxiliar aislados eq2 kurosh en bola}
\{\zeta_1,\dots,\zeta_r, \alpha_1,\dots \alpha_s, \beta_1,\dots,\beta_t\}\subseteq \bola_{\equisg\cup \equish}(n),
\end{equation}
\begin{equation*}
\ga \lambda_{n+1} \ce\neq\lambda_{n+1} \ce\;\mbox{ , }\;\gb \lambda_{n+1} \ce\neq\lambda_{n+1} \ce\;\mbox{ and }\;\ga^{\varepsilon_\gr} \lambda_{n+1} \ce\neq \gb^{\varepsilon_\grb}\lambda_{n+1} \ce\;\;\mbox{ for every }\; \varepsilon_\gr,\varepsilon_\grb\in\{\pm 1\}.
\end{equation*}
From this last fact, and by changing $\ga$ or $\gb$ (or both) by its inverse if necessary, we can assume that $\lambda_{n+1}^+\ce\lle \ga\lambda_{n+1}^+ \ce$, that $\lambda_{n+1}^+\ce\lle \gb\lambda_{n+1}^+ \ce$ and that $\ga \lambda_{n+1} \ce\neq \gb\lambda_{n+1} \ce$. Notice that, by interchanging the groups $\gr$ and $\grb$ by symmetry if necessary, we can also suppose that $\gb\lambda_{n+1}^+\ce\lle \ga\lambda_{n+1}^+ \ce$.
Hence, we have 
\begin{equation}\label{teorema auxiliar aislados eq3 lambdas inecuacion}
\lambda_{n+1}^+\ce\lle \gb\lambda_{n+1}^+ \lle \ga\lambda_{n+1}^+\ce
\end{equation}
We now define $$\gamma_*=(\gb\lambda_{n+1}^+)\inv \ga\lambda_{n+1}^+,$$ so we have $$\ce\lle \gamma_* \ce.$$ Notice that $\ga,\gb,\lambda_{n+1}^+\in \gen{\equisg,\equish}$, so we also have $\gamma_*\in \gen{\equisg,\equish}$.\\

At this point, we have built two finite sets $\equisg\subseteq \gr$ and $\equish\subseteq \grb$ and an element $\gamma_{*}\in \gen{\equisg,\equish}$ such that $\ce\lle \gamma_{*}\ce$. In order to conclude the proof, we are going to prove that the triple $(\equisg, \equish, \gamma_{*})$ satisfies the property claimed in the statement. In other words, we need to prove that for all finitely generated subgroups $\gr_1\subseteq \gr$ and $\grb_1\subseteq \grb$ containing $\equisg$ and $\equish$ respectively, the group $\Gamma=\gr_1*\grb_1$ admits a left-preorder $\lleq^*$ relative to $\ce_\Gamma=\Gamma\cap \ce$ such that $\gamma_{*}\ce_\Gamma\lle^* \ce_\Gamma$ and $\ce_\Gamma\lle^* \elefiny\ce_\Gamma$ for all $\elefiny\in\conjfiny$.\\

Consider $\gr_1\subseteq \gr$ and $\grb_1\subseteq \grb$  two finitely generated subgroups containing $\equisg$ and $\equish$ respectively. Consider $\Gamma=\gr_1*\grb_1$, which is finitely generated. Define $\ce_\Gamma=\Gamma\cap \ce$. We are going to build the left-preorder $\lleq^*$ by using the point of view of dynamical realizations and then applying \Cref{extension desde dinamica} to obtain the needed left-preorder.\\

The restriction of $\lleq$ to $\Gamma$ is a left-preorder on $\Gamma$ relative to $\ce_\Gamma=\ce\cap \Gamma$, which we will denote by $\lleq$ as an abuse of notation. Notice that this left-preorder is given by $\gamma_1\ce\lleq \gamma_2\ce$ if and only if $\gamma_1\ce_\Gamma\lleq \gamma_2\ce_\Gamma$ for all $\gamma_1,\gamma_2\in \Gamma$. Since $\Gamma$ is finitely generated, by \Cref{lema realizacion dinamica orden}, we can consider $\df{D}{\Gamma}{\homeomr}$ a dynamical realization-like homomorphism for $\lleq$ as left-preorder on $\Gamma$ relative to $\ce_\Gamma$ with reference point $0$. Then, we have 
\begin{equation}\label{teorema auxiliar aislados eq4 D es realizacion dinamica}
\gamma_1\ce_\Gamma\lleq \gamma_2\ce_\Gamma \Leftrightarrow D(\gamma_1)(0)\leq D(\gamma_2)(0) \;\mbox{ for all }\; \gamma_1,\gamma_2\in \Gamma.
\end{equation}
From \Cref{teorema auxiliar aislados eq3 lambdas inecuacion}, we deduce that $$D(\lambda_{n+1}^+)(0) < D(\gb\lambda_{n+1}^+)(0) < D(\ga\lambda_{n+1}^+)(0).$$ By density of $\real$, we can consider $\puntox_0,\puntox_1,\puntoy_0,\puntoy_1\in\real$ such that
\begin{equation}\label{teorema auxiliar aislados eq5 colocacion de los puntos}
D(\lambda_{n+1}^+)(0) < \puntox_0 < \puntox_1 < D(\gb\lambda_{n+1}^+)(0) < D(\ga\lambda_{n+1}^+)(0)< \puntoy_1 < \puntoy_0.
\end{equation}
The support of a map $\df{\mu}{\real}{\real}$ is defined to be the set $\{\puntox\in \real\tq \mu(\puntox)\neq \puntox\}$. We take $\phi\in \homeomr$ such that its support is equal to $(\puntox_0,\puntoy_0)$ and satisfying that $\phi(\puntox_1)>\puntoy_1$. We claim that
\begin{equation}\label{teorema auxiliar aislados eq6 phi de 0}
\phi(0)=0\;\mbox{ and }\;\phi\inv(0)=0.
\end{equation}
Indeed, since $\ce_\Gamma\lle\lambda_{n+1}^+\ce_\Gamma$, then $0<D(\lambda_{n+1}^+)(0)$, and it implies that $0\notin (\puntox_0,\puntoy_0)$. In particular, $0$ is not on the support of $\phi$, so this proves \Cref{teorema auxiliar aislados eq6 phi de 0}.\\

Recall that $\Gamma=\gr_1*\grb_1$. By the universal property of the free product, we can build $\df{D_\phi}{\Gamma}{\homeomr}$ a group homomorphism defined by $D_\phi(\overline{\ga})=D(\overline{\ga})$ for all $\overline{\ga}\in \gr_1$ and $D_\phi(\overline{\gb})=\phi\co D(\overline{\gb})\co\phi\inv$ for all $\overline{\gb}\in \grb_1$.\\

To make the proof easier to follow, we are going to conclude the proof by assuming three claims, and we are going to prove them when we are done with the main proof. The three equations we are going to assume are  
\begin{equation}\label{teorema auxiliar aislados eq a probar 1 coincide en la bola}
D_\phi(w)(0)=D(w)(0)\; \mbox{ for all } w\in \bola_{\equisg\cup \equish}(n+1)\mbox{,}
\end{equation}
\begin{equation}\label{teorema auxiliar aislados eq a probar 2 inecuacion con lambda}
D_\phi(\ga\lambda_{n+1}^+)(0)< D_\phi(\gb\lambda_{n+1}^+)(0)  \;\mbox{ and }
\end{equation}
\begin{equation}\label{teorema auxiliar aislados eq a probar 3 contenido estabilizador}
\ce_\Gamma\subseteq \stab_{D_\phi}(0).
\end{equation}
We are able to build the needed left-preorder using $D_\phi$. Combining \Cref{teorema auxiliar aislados eq a probar 3 contenido estabilizador} and \Cref{extension desde dinamica}, we obtain a left-preorder $\lleq^*$ on $\Gamma$ relative to $\ce_\Gamma$ such that 
\begin{equation}\label{teorema auxiliar aislados eq final}
D_\phi(\gamma_1)(0)<D_\phi(\gamma_2)(0) \;\Longrightarrow \; \gamma_1\ce_\Gamma\lle^* \gamma_2\ce_\Gamma \quad \mbox{for all} \gamma_1,\gamma_2\in\Gamma.
\end{equation}
To conclude the proof, we only need to prove that $\gamma_{*}\ce_\Gamma\lle^* \ce_\Gamma$ and $\ce_\Gamma\lle^* \elefiny\ce_\Gamma$ for all $\elefiny\in\conjfiny$. We check each condition. From \Cref{teorema auxiliar aislados eq a probar 2 inecuacion con lambda} and \Cref{teorema auxiliar aislados eq final} we deduce that $$\ga\lambda_{n+1}^+\ce_\Gamma\lle^* \gb\lambda_{n+1}^+\ce_\Gamma.$$ Then, we deduce $$\ce_\Gamma\lle^* (\ga\lambda_{n+1}^+)\inv \gb\lambda_{n+1}^+\ce_\Gamma=(\gamma_*)\inv\ce_\Gamma.$$ This proves the condition $\gamma_{*}\ce_\Gamma\lle^* \ce_\Gamma$. Now, given $\elefiny\in\conjfiny$, by hypothesis we know that $\ce\lle \elefiny\ce$. Since $D$ is the dynamical realization-like homomorphism for the restriction of $\lleq$ with reference point $0$, we deduce that $$0<D(\elefiny)(0).$$ Combining this fact with \Cref{teorema auxiliar aislados eq1 x en bola} and \Cref{teorema auxiliar aislados eq a probar 1 coincide en la bola}, we deduce $$0<D_\phi(\elefiny)(0).$$ Applying \Cref{teorema auxiliar aislados eq final}, we obtain $$\ce_\Gamma\lle^* \elefiny\ce_\Gamma.$$ This proves the condition $\ce_\Gamma\lle^* \elefiny\ce_\Gamma$ for all $\elefiny\in\conjfiny$. Therefore, we have shown $\gamma_{*}\ce_\Gamma\lle^* \ce_\Gamma$ and $\ce_\Gamma\lle^* \elefiny\ce_\Gamma$ for all $\elefiny\in\conjfiny$. This finishes the main proof of the theorem.\\

Recall that, for the proof, we have assumed that \Cref{teorema auxiliar aislados eq a probar 1 coincide en la bola}, \Cref{teorema auxiliar aislados eq a probar 2 inecuacion con lambda} and \Cref{teorema auxiliar aislados eq a probar 3 contenido estabilizador} hold. Now, we are going to prove each statement.\\ 

We start by proving \Cref{teorema auxiliar aislados eq a probar 1 coincide en la bola}. We claim that
\begin{equation}\label{teorema auxiliar aislados primera prueba eq 1}
D_\phi(\overline{\gb})(\puntox)=D(\overline{\gb})(\puntox)\; \mbox{ for all } \overline{\gb}\in \equish \;\mbox{ and all } \puntox\leq D(\lambda_{n}^+)(0).
\end{equation}
Indeed, consider $\overline{\gb}\in \equish$ and $\puntox\leq D(\lambda_{n}^+)(0)$. Notice that have $\overline{\gb}\lambda_{n}^+\in \bola_{\equisg\cup \equish}(n+1)$. Then, by definition of $\lambda_{n}^+$, we have $$\overline{\gb}\lambda_{n}^+\ce_\Gamma\lleq \lambda_{n+1}^+\ce_\Gamma.$$ Then, using \Cref{teorema auxiliar aislados eq4 D es realizacion dinamica}, we obtain $$D(\overline{\gb}\lambda_{n}^+)(0)\leq D(\lambda_{n+1}^+)(0).$$ Combining this with $\puntox\leq D(\lambda_{n}^+)(0)$, with \Cref{teorema auxiliar aislados eq5 colocacion de los puntos} and with the fact that $D(\overline{\gb})$ is order preserving, we obtain $$D(\overline{\gb})(\puntox)\leq D(\overline{\gb}) \left(D(\lambda_{n}^+)(0)\right)=D(\overline{\gb}\lambda_{n}^+)(0)\leq D(\lambda_{n+1}^+)(0)<\puntox_0.$$ Hence, $D(\overline{\gb})(\puntox)<\puntox_0$. Notice also that $\puntox<\puntox_0$ as a consequence of $\puntox\leq D(\lambda_{n}^+)(0)$ and \Cref{teorema auxiliar aislados eq5 colocacion de los puntos}. We have shown that $D(\overline{\gb})(\puntox)$ and $\puntox$ are both smaller than $\puntox_0$, so they are not at the support of $\phi$. Using this and the definition of $D_\phi$, we have $ D_\phi(\overline{\gb})(\puntox)=D(\overline{\gb})(\puntox)$. This proves \Cref{teorema auxiliar aislados primera prueba eq 1}, as we claimed.\\

By definition of $D_\phi$, we have $$D_\phi(\overline{\ga})(\puntox)=D(\overline{\ga})(\puntox)\; \mbox{ for all } \overline{\ga}\in \equisg \;\mbox{ and all } \puntox\leq D(\lambda_{n}^+)(0).$$ Combining this with \Cref{teorema auxiliar aislados primera prueba eq 1}, we deduce that, for all $\gamma\in \equisg\cup \equish$, we have that $D(\gamma)$ and $D_\phi(\gamma)$ are equal as maps when restricted to $[D(\lambda_n^-)(0),D(\lambda_n^+)(0)]$. Then, from this fact and from \Cref{lema realizaciones}, we deduce \Cref{teorema auxiliar aislados eq a probar 1 coincide en la bola}, as we wanted. \\

We are now going to prove \Cref{teorema auxiliar aislados eq a probar 2 inecuacion con lambda}. Recall that $g$ and $h$ are the fixed non-trivial elements of $\equisg$ and $\equish$. From \Cref{teorema auxiliar aislados eq5 colocacion de los puntos} and \Cref{teorema auxiliar aislados eq6 phi de 0}, we obtain $$\left(D(\gb\lambda_{n+1}^+)\co \phi\inv\right)(0)=D(\gb\lambda_{n+1}^+)\left(\phi\inv(0)\right)=D(\gb\lambda_{n+1}^+)(0)>\puntox_1.$$ Applying that $\phi$ is order preserving and $\phi(\puntox_1)>\puntoy_1$ by definition, we deduce $$\left(\phi\co D(\gb\lambda_{n+1}^+)\co \phi\inv\right)(0)>\phi(\puntox_1)>\puntoy_1.$$ Combining this with \Cref{teorema auxiliar aislados eq5 colocacion de los puntos}, we deduce $$D(\ga\lambda_{n+1}^+)(0)<\puntoy_1<\left(\phi\co D(\gb\lambda_{n+1}^+)\co \phi\inv\right)(0).$$
Hence, we have shown that 
\begin{equation}\label{teorema auxiliar aislados segunda prueba eq 1}
D(\ga\lambda_{n+1}^+)(0)<\left(\phi\co D(\gb\lambda_{n+1}^+)\co \phi\inv\right)(0).
\end{equation}
Now since $\lambda_{n+1}^+\in \bola_{\equisg\cup \equish}(n+1)$, applying \Cref{teorema auxiliar aislados eq a probar 1 coincide en la bola}, we obtain 
\begin{equation}\label{teorema auxiliar aislados segunda prueba eq 2}
D_\phi(\lambda_{n+1}^+)(0)=D(\lambda_{n+1}^+)(0).
\end{equation}
By definition of $D_\phi$ knowing that $\ga\in \gr_1$, we obtain $D_\phi(\ga)=D(\ga)$. From this and from \Cref{teorema auxiliar aislados segunda prueba eq 2} we deduce $$D_\phi(\ga\lambda_{n+1}^+)(0)=D_\phi(\ga)\left(D_\phi(\lambda_{n+1}^+)(0)\right)=D(\ga)\left(D(\lambda_{n+1}^+)(0)\right)=D(\ga\lambda_{n+1}^+)(0).$$ Then, we have shown that 
\begin{equation}\label{teorema auxiliar aislados segunda prueba eq 3}
D_\phi(\ga\lambda_{n+1}^+)(0)=D(\ga\lambda_{n+1}^+)(0).
\end{equation}
From \Cref{teorema auxiliar aislados eq5 colocacion de los puntos} we notice that $D(\lambda_{n+1}^+)(0)$ is not in the support of $\phi$, so then
\begin{equation}\label{teorema auxiliar aislados segunda prueba eq 4}
\phi\inv\left(D(\lambda_{n+1}^+)(0)\right)=D(\lambda_{n+1}^+)(0).
\end{equation}
By definition of $D_\phi$ knowing that $\gb\in \grb_1$, we obtain $$D_\phi(\gb)=(\phi \co D\co \phi\inv)(\gb).$$ From this, \Cref{teorema auxiliar aislados eq6 phi de 0}, \Cref{teorema auxiliar aislados segunda prueba eq 4} and \Cref{teorema auxiliar aislados segunda prueba eq 2} we deduce \begin{equation*}
\begin{split} 
D_\phi(\gb\lambda_{n+1}^+)(0)=D_\phi(\gb)\left(D_\phi(\lambda_{n+1}^+)(0)\right)=D_\phi(\gb)\left(D(\lambda_{n+1}^+)(0)\right)=\left(\phi\co D(\gb)\co \phi\inv\right)\left(D(\lambda_{n+1}^+)(0)\right)=\\
=\left(\phi\co D(\gb)\right)\left(D(\lambda_{n+1}^+)(0)\right)=\left(\phi\co D(\gb\lambda_{n+1}^+)\right)(0)=\left(\phi \co D(\gb\lambda_{n+1}^+)\right)(\phi\inv (0))=\left(\phi \co D(\gb\lambda_{n+1}^+)\co \phi\inv\right) (0).
\end{split}
\end{equation*} Then, we have proven that 
\begin{equation}\label{teorema auxiliar aislados segunda prueba eq 5}
D_\phi(\gb\lambda_{n+1}^+)(0)=\left(\phi \co D(\gb\lambda_{n+1}^+)\co \phi\inv\right) (0).
\end{equation}
We can combine \Cref{teorema auxiliar aislados segunda prueba eq 3}, \Cref{teorema auxiliar aislados segunda prueba eq 1} and \Cref{teorema auxiliar aislados segunda prueba eq 5}, so we obtain $$D_\phi(\ga\lambda_{n+1}^+)(0)=D(\ga\lambda_{n+1}^+)(0)<\left(\phi \co D(\gb\lambda_{n+1}^+)\co \phi\inv \right)(0)=D_\phi(\gb\lambda_{n+1}^+)(0).$$
This proves \Cref{teorema auxiliar aislados eq a probar 2 inecuacion con lambda}, as we wanted.\\

We finally prove \Cref{teorema auxiliar aislados eq a probar 3 contenido estabilizador}. In order to prove it, recall that, by definition, we have that $\fino\subseteq \gen{\equisg}*\gen{\equish}$ and that $\gr_1\subseteq \gr$ and $\grb_1\subseteq \grb$ are subgroups satisfying $\equisg\subseteq \gr_1$ and $\equish\subseteq \grb_1$. In particular, $\fino$ is contained in $\gen{\gr_1,\grb_1}$. Then, applying Fact I, we deduce that $$\ce_\Gamma=\gen{\{\zeta_i \tq i=1,\dots,r\}\cup \left(\bigcup_{j=1}^s\alpha_j\grfaca_j\alpha_j\inv \right)\cup \left(\bigcup_{l=1}^t\beta_l\grfacb_l\beta_l\inv \right)}$$ where $\grfaca_j$ is  subgroup of $\gr_0$ for $j=1,\dots,s$ and where $\grfacb_t$ is a subgroup of $\grb_0$ for $t=1,\dots,l$. From this description of $\ce_\Gamma$, in order to prove \Cref{teorema auxiliar aislados eq a probar 3 contenido estabilizador}, we only need to prove the following three equations involving the generators of $\ce_\Gamma$:
\begin{equation}\label{teorema auxiliar aislados tercera prueba eq 1}
\zeta_i\in \stab_{D_\phi}(0) \;\mbox{ for all } i=1,\dots,r \;,
\end{equation}
\begin{equation}\label{teorema auxiliar aislados tercera prueba eq 2}
\alpha_ja\alpha_j\inv\in \stab_{D_\phi}(0) \;\mbox{ for all } j=1,\dots,s \;\mbox{ and all }a\in\grfaca_j\;,
\end{equation}
\begin{equation}\label{teorema auxiliar aislados tercera prueba eq 3}
\beta_lb\beta_l\inv\in \stab_{D_\phi}(0) \;\mbox{ for all } l=1,\dots,t \;\mbox{ and all }b\in\grfacb_l.
\end{equation}
We are going to prove each equation separately. Before doing that, recall that $$\ce_\Gamma=\stab_{D}(0)$$ since $D$ is a dynamical realization-like homomorphism for a left-preorder relative to $\ce_\Gamma$ with reference point $0$.\\

We are going to prove \Cref{teorema auxiliar aislados tercera prueba eq 1}. Take $i\in\{1,\dots,r\}$. Since $\ce_\Gamma=\stab_{D}(0)$ and $\zeta_i\in \ce_\Gamma$, we obtain $D(\zeta_i)(0)=0$. The previous equation combined with \Cref{teorema auxiliar aislados eq2 kurosh en bola} and \Cref{teorema auxiliar aislados eq a probar 1 coincide en la bola} give us that $D_\phi(\zeta_i)(0)=D(\zeta_i)(0)$ is equal to $0$. This shows \Cref{teorema auxiliar aislados tercera prueba eq 1}, as we wanted.\\

We are going to show \Cref{teorema auxiliar aislados tercera prueba eq 2}. Take $j\in\{1,\dots,s\}$ and $\gfa\in\grfaca_j$. In order to prove \Cref{teorema auxiliar aislados tercera prueba eq 2} we are going to show $$D_\phi(\alpha_j\inv)(0)=D_\phi(\gfa\alpha_j\inv)(0),$$ as they are equivalent statements because $D_\phi$ is an homomorphism. On the one hand, \Cref{teorema auxiliar aislados eq2 kurosh en bola} implies
\begin{equation*}
\alpha_j\inv\in \bola_{\equisg\cup \equish}(n).
\end{equation*}
Combining this with \Cref{teorema auxiliar aislados eq a probar 1 coincide en la bola}, we deduce
\begin{equation}\label{teorema auxiliar aislados tercera prueba eq 2.1}
D_\phi(\alpha_j\inv)(0)=D(\alpha_j\inv)(0).
\end{equation}
On the other hand, as $\alpha_j\gfa\alpha_j\inv\in \ce_\Gamma$ and $\ce_\Gamma=\stab_{D}(0)$, using that $D$ is an homomorphism we deduce
\begin{equation}\label{teorema auxiliar aislados tercera prueba eq 2.2}
D(\alpha_j\inv)(0)=D(\gfa\alpha_j\inv)(0).
\end{equation}
Also, as $\gfa\in \gr_1$, from the definition of $D_\phi$ we obtain 
\begin{equation}\label{teorema auxiliar aislados tercera prueba eq 2.3}
D_\phi(\gfa)=D(\gfa).
\end{equation}
Combining \Cref{teorema auxiliar aislados tercera prueba eq 2.1}, \Cref{teorema auxiliar aislados tercera prueba eq 2.2} and \Cref{teorema auxiliar aislados tercera prueba eq 2.3}, and using that $D$ and $D_\phi$ are homomorphisms, we deduce that
\begin{equation*}
D_\phi(\gfa\alpha_j\inv)(0)=D_\phi(\gfa)\left(D_\phi(\alpha_j\inv)(0)\right)=D(\gfa)\left(D(\alpha_j\inv)(0)\right)=D(\gfa\alpha_j\inv)(0)=D(\alpha_j\inv)(0)=D_\phi(\alpha_j\inv)(0)
\end{equation*}
Hence, we have proven $D_\phi(\alpha_j\inv)(0)=D_\phi(\gfa\alpha_j\inv)(0)$, as we wanted. Therefore, we have shown \Cref{teorema auxiliar aislados tercera prueba eq 2}.\\

We are going to prove \Cref{teorema auxiliar aislados tercera prueba eq 3}. Consider $l\in\{1,\dots,t\}$ and $\gfb\in\grfacb_l$. Analogously as before, in order to prove \Cref{teorema auxiliar aislados tercera prueba eq 3} we are going to prove $$D_\phi(\beta_l\inv)(0)=D_\phi(\gfb\beta_l\inv)(0),$$ as both are equivalent statements as $D_\phi$ is an homomorphism. Previously, we are going to derive some useful equations that, when combined, will let us obtain the expression we needed. From \Cref{teorema auxiliar aislados eq2 kurosh en bola}, we obtain
\begin{equation}\label{teorema auxiliar aislados tercera prueba eq 3.1}
\beta_l\inv\in \bola_{\equisg\cup \equish}(n).
\end{equation}
This implies, by \Cref{teorema auxiliar aislados eq a probar 1 coincide en la bola}, that
\begin{equation}\label{teorema auxiliar aislados tercera prueba eq 3.2}
D_\phi(\beta_l\inv)(0)=D(\beta_l\inv)(0).
\end{equation}
Now, from \Cref{teorema auxiliar aislados tercera prueba eq 3.1} and the definition of $\lambda^+_{n+1}$ we obtain $$\beta_l\inv \ce_\Gamma \lle \lambda_{n+1}^+ \ce_\Gamma.$$ Combining this with \Cref{teorema auxiliar aislados eq4 D es realizacion dinamica} and \Cref{teorema auxiliar aislados eq5 colocacion de los puntos}, we obtain $$D(\beta_l\inv)(0)<D(\lambda_{n+1}^+)(0) < \puntox_0.$$ In particular, $D(\beta_l\inv)(0)$ is not in the support of $\phi$, so then 
\begin{equation}\label{teorema auxiliar aislados tercera prueba eq 3.3}
\phi\inv\left(D(\beta_l\inv)(0)\right)=D(\beta_l\inv)(0)\;\mbox{ and }\; \phi\left(D(\beta_l\inv)(0)\right)=D(\beta_l\inv)(0).
\end{equation}
On the other hand, as we know that $\beta_l\gfb\beta_l\inv\in \ce_\Gamma$ and $\ce_\Gamma=\stab_{D}(0)$, using that $D$ is an homomorphism we obtain
\begin{equation}\label{teorema auxiliar aislados tercera prueba eq 3.4}
D(\beta_l\inv)(0)=D(\gfb\beta_l\inv)(0).
\end{equation}
Also, as $\gfb\in \grb_1$, from the definition of $D_\phi$ we obtain 
\begin{equation}\label{teorema auxiliar aislados tercera prueba eq 3.5}
D_\phi(\gfb)=\phi\co D(\gfb)\co\phi\inv.
\end{equation}
Combining \Cref{teorema auxiliar aislados tercera prueba eq 3.2}, \Cref{teorema auxiliar aislados tercera prueba eq 3.3}, \Cref{teorema auxiliar aislados tercera prueba eq 3.4} and \Cref{teorema auxiliar aislados tercera prueba eq 3.5}, and using that $D$ and $D_\phi$ are homomorphisms, we obtain the chain of equalities
\begin{equation*}
\begin{split} 
D_\phi(\gfb\beta_l\inv)(0)=D_\phi(\gfb)\left(D_\phi(\beta_l\inv)(0)\right)=\left( \phi\co D(\gfb)\phi\inv\right) \left(D(\beta_l\inv)(0)\right)=\left( \phi\co D(\gfb)\right) \left(\phi\inv \left(D(\beta_l\inv)(0)\right)\right)=\\
=\left( \phi\co D(\gfb)\right) \left( D(\beta_l\inv)(0)\right)=\phi \left( D(\gfb\beta_l\inv)(0) \right)=\phi \left( D(\beta_l\inv)(0) \right)=D(\beta_l\inv)(0)=D_\phi(\beta_l\inv)(0)
\end{split} 
\end{equation*}
This proves that $D_\phi(\beta_l\inv)(0)=D_\phi(\gfb\beta_l\inv)(0)$, as we wanted. Therefore, we have shown \Cref{teorema auxiliar aislados tercera prueba eq 3}.\\
This concludes the proof of \Cref{teorema auxiliar aislados eq a probar 3 contenido estabilizador}, so we can finish the proof.
\end{proof}

Next corollary follows by applying Theorem 27.7 of \cite{topology-munkres}, knowing that $\pord_\ce(\gr*\grb)$ is compact, Hausdorff and has no isolated points.

\begin{coro}
Let $\gr$ and $\grb$ be two non-trivial groups. Let $\ce$ be a proper subgroup of $\gr*\grb$ having finite Kurosh rank. If $\ce$ is left-relatively convex in $\gr*\grb$, then $\pord_\ce(\gr*\grb)$ is uncountable.
\end{coro}

Next result follows by applying the characterization of Cantor sets as a compact, Hausdorff, totally disconnected, non-empty metric space having no isolated points (\Cref{cuando es cantor}).

\begin{coro}\label{espacio de ordenes es un conjunto de Cantor}
Let $\gr$ and $\grb$ be two finitely generated non-trivial groups. Let $\ce$ be a proper subgroup of $\gr*\grb$ having finite Kurosh rank. If $\ce$ is left-relatively convex in $\gr*\grb$, then $\pord_\ce(\gr*\grb)$ is homeomorphic to the Cantor set.
\end{coro}

Last corollary is also true if we change the condition of $\gr$ and $\grb$ being finitely generated by the condition of $\gr$ and $\grb$ being countable. In particular, it is also true if we change that condition by $\gr$ and $\grb$ being countable.\\

We can state now the theorem about $\pord(\gr*\grb)$.

\begin{teo}\label{teorema aislados general}
Let $\gr$ and $\grb$ be two non-trivial groups. Then, $\pord(\gr*\grb)$ has no isolated elements.
\begin{proof}
We want to show that, if $\lleq\in\pord(\gr*\grb)$, then $\lleq$ is not isolated in $\pord(\gr*\grb)$. Indeed, consider $\lleq$ an element of $\pord(\gr*\grb)$. Then, $\lleq$ is a left-preorder on $\gr*\grb$ relative to $\ce$ for some proper subgroup $\ce$ of $\gr*\grb$, that is, $\lleq$ is an element of the subspace $\pord_\ce(\gr*\grb)$. We are going to prove that $\lleq$ is not isolated $\pord(\gr*\grb)$ by considering two cases: when the Kurosh rank of $\ce$ with respect to $\gr*\grb$ is finite or infinite. If $\ce$ has finite Kurosh rank with respect to $\gr*\grb$, by \Cref{teorema aislados}, we know that $\lleq$ is not isolated on $\pord_\ce(\gr*\grb)$, so $\lleq$ is not isolated on $\pord(\gr*\grb)$. Hence, we can assume that we are in the case where $\ce$ has infinite Kurosh rank with respect to $\gr*\grb$.\\

Consider $\conjfinf\subseteq \gr*\grb$ a finite subset. The fact that $\lleq$ is not isolated on $\pord(\gr*\grb)$ will follow if we build a left-preorder $\lleq'$ that is non equal to $\lleq$ but both coinciding on the set $\conjfinf$. Consider $P$ the positive cone of $\lleq$. The set $\conjfinf$ is finite and $\ce$ has infinite Kurosh rank with respect to $\gr*\grb$, so there is a proper free factor $\ce_1$ of $\ce$ such that $\conjfinf \cap \ce$ is contained in $\ce_1$ (\Cref{obs bs rango kurosh infinito}). Then, we can write $$\ce = \ce_1 * \ce_2.$$ Since the subgroup $\ce_1$ is left-relatively convex on $\ce$ (\Cref{ordenar con factor convexo}), we deduce from \Cref{extension orden} that there is $\lleq'$ a left-preorder on $\gr*\grb$ relative to $\ce_1$ such that $\gx_1 \ce\lle \gx_2 \ce$ implies $\gx_1\ce_1\lle' \gx_2\ce_1$. In particular, $\lleq$ and $\lleq'$ coincide on $\conjfinf\settminus \ce$. Notice that $\lleq$ and $\lleq'$ coincide on $\conjfinf\cap \ce$ since $\lleq'$ is a left-preorder on $\gr*\grb$ relative to $\ce_1$ and $\conjfinf \cap \ce$ is contained in $\ce_1$. Therefore, $\lleq$ and $\lleq'$ coincide on $\conjfinf$. The fact that they are not equal is direct since $\ce_1$ and $\ce$ are not equal and since $\lleq$ is a left-preorder on $\gr*\grb$ relative to $\ce$ and $\lleq'$ is a left-preorder on $\gr*\grb$ relative to $\ce_1$.
\end{proof}
\end{teo}

We have some consequences of the theorem.

\begin{coro}\label{corolario espacio de preordenes es incontable producto libre}
Let $\gr$ and $\grb$ be two non-trivial groups. If $\pord(\gr*\grb)$ is non-empty, then it is uncountable.
\begin{proof}
Assume that $\pord(\gr*\grb)$ is non-empty. Let $\lleq_{T}$ be the trivial preorder. We are going to consider two cases depending on if $\lleq_{T}$ is isolated on the set $\pord (\gr*\grb)\cup\{\lleq_{T}\}$ or not.\\

Assume that we are on the case where $\lleq_{T}$ is not isolated on the set $\pord (\gr*\grb)\cup\{\lleq_{T}\}$. Then, by \Cref{teorema aislados general}, we deduce that $\pord (\gr*\grb)\cup\{\lleq_{T}\}$ has no isolated elements, and by \Cref{espacio de ordenes con el trivial es cerrado} and \Cref{espacio de ordenes es totalmente disconexo y hausdorff}, we deduce that $\pord (\gr*\grb)\cup\{\lleq_{T}\}$ is Hausdorff and has no isolated points. Applying Theorem 27.7 of \cite{topology-munkres}, we deduce that $\pord (\gr*\grb)\cup\{\lleq_{T}\}$ is uncountable. Therefore, $\pord(\gr*\grb)$ is uncountable.\\

Assume that we are on the case where $\lleq_{T}$ is isolated on the set $\pord (\gr*\grb)\cup\{\lleq_{T}\}$. This means that $\pord (\gr*\grb)$ is closed on $\pord (\gr*\grb)\cup\{\lleq_{T}\}$, but $\pord (\gr*\grb)\cup\{\lleq_{T}\}$ is compact by \Cref{espacio de ordenes con el trivial es cerrado}, so then $\pord (\gr*\grb)$ is compact. Also, by \Cref{espacio de ordenes es totalmente disconexo y hausdorff} and \Cref{teorema aislados general}, we deduce that $\pord (\gr*\grb)$ is Hausdorff and has no isolated points. Then, by Theorem 27.7 of \cite{topology-munkres}, we deduce that $\pord(\gr*\grb)$ is uncountable.
\end{proof}
\end{coro}

Next corollary follows by applying the characterization of Cantor sets (\Cref{cuando es cantor}) and by knowing that the free product of left-orderable groups is left-orderable (see Corollary 6.1.3 of \cite{kopytovmedevev}).

\begin{coro}\label{corolario espacio de preordenes es cantor producto libre}
Let $\gr$ and $\grb$ be two non-trivial finitely generated left-orderable groups. Then, $\pord(\gr*\grb)$ is homeomorphic to the Cantor set.
\end{coro}

From the theorem we obtain a description of this spaces for non-abelian free groups. Recall that, for a subgroup, being finitely generated implies having a finite Kurosh rank.

\begin{coro}\label{descripcion grupo libre}
Let $\mathbb{F}$ be a non-abelian free group. Then,
\begin{enumerate}
\item $\pord_\ce(\mathbb{F})$ has no isolated elements for each $\ce$ finitely generated left-relatively convex subgroup on $\mathbb{F}$. In particular, if $\mathbb{F}$ has at most countable free rank, then $\pord_\ce(\mathbb{F})$ is a Cantor set.
\item $\pord(\mathbb{F})$ has no isolated elements. In particular, if $\mathbb{F}$ has finite free rank, then $\pord(\mathbb{F})$ is a Cantor set.
\end{enumerate}
\end{coro}

As a corollary of the theorem, when combined with \Cref{fg como semigrupo implica asilado}, we obtain the following property about positive cones.

\begin{coro}\label{cono positivo no finitamente generado como corolario de aislado}
Let $\gr$ and $\grb$ be two non-trivial groups. Let $\ce$ be a left-relatively convex subgroup on $\gr*\grb$ having finite Kurosh rank with respect to $\gr*\grb$. Let $\lleq$ be a left-preorder on $\gr*\grb$ relative to $\ce$. Then, the positive cone of $\lleq$ is not finitely generated as a subsemigroup of $\gr*\grb$.
\end{coro}

\section{Auxiliary reults on subgroups of finite Kurosh rank}

Recall that, in the proof of \Cref{teorema auxiliar aislados}, we used two facts when assuming that $\gr$ and $\grb$ are two non-trivial groups and $\ce$ is a subgroup of  $\gr*\grb$ having finite Kurosh rank with respect to $\gr*\grb$. These facts are:
\begin{itemize}
\item[\textbf{Fact I:}] There exists $\fino,\{\zeta_i\}_{i=1}^r,\{\alpha_j\}_{j=1}^s,\{\beta_t\}_{t=1}^l\subseteq \gr*\grb$ finite sets satisfying the following property: for every $\gr_0$ subgroup of $\gr$ and every $\grb_0$ subgroup of $\grb$ satisfying $\fino\subseteq \gen{\gr_0,\grb_0}$ there exists $\grfaca_j$ subgroup of $\gr_0$ for each $j=1,\dots,s$ and there exists $\grfacb_t$ subgroup of $\grb_0$ for each $t=1,\dots,l$ such that 
\begin{equation*}
\ce\cap\gen{\gr_0,\grb_0}=\gen{\{\zeta_i \tq i=1,\dots,r\}\cup \left(\bigcup_{j=1}^s\alpha_j\grfaca_j\alpha_j\inv \right)\cup \left(\bigcup_{l=1}^t\beta_l\grfacb_l\beta_l\inv \right)}.
\end{equation*}
\item[\textbf{Fact II:}] Let $(\lambda_l)_{l\in\nat}$ be a sequence of elements in $\gr*\grb$ such that $\lambda_l \ce\neq\lambda_m \ce$ for all $l\neq m$. Then, there exists $n\in\nat$ such that for all $k\geq n$ there are $\ga\in \equisg$ and $\gb\in \equish$ satisfying $\ga^{\varepsilon_\gr} \lambda_k \ce\neq\lambda_k \ce$, $\gb^{\varepsilon_\grb} \lambda_k \ce\neq\lambda_k \ce$ and $\ga^{\varepsilon_\gr} \lambda_k \ce\neq \gb^{\varepsilon_\grb}\lambda_k \ce$ for every $\varepsilon_\gr,\varepsilon_\grb\in\{\pm 1\}$.
\end{itemize}

In this section we use Bass-Serre theory to prove these facts. Our basic reference for Bass-Serre theory is \cite{serre}. The section is a bit technical. For a non-deep lecture, one can skip this section by just assuming the previous facts and by just considering the Kurosh rank as a black box that satisfies that properties. Fact I is shown in \Cref{bass serre para interseccion} and Fact II is shown in \Cref{corolario final lambda g y lambda h}. However, it is remarkable that \Cref{bass serre para interseccion} could be an interesting result by itself. As far as the author knows, that result has not been considered in the literature, and it provides a description for certain intersections on free products and, as a consequence, one obtain a bound for the Kurosh rank of that intersection.\\

We start by fixing basic notation from graph theory. A \textit{graph} $\grf$ consists of a non-empty set of \textit{vertices} $\grfv$, a set of \textit{edges} $\grfe$, an \textit{inverse-edge map} $\df{\inv}{\grfe}{\grfe}$ which is an involution without fixed-points, a \textit{initial-vertex map} $\df{\origen}{\grfe}{\grfv}$ and a \textit{terminal-vertex map} $\df{\final}{\grfe}{\grfv}$ satisfying $\origen(\eda)=\final(\eda\inv)$ for all $\eda\in\grfe$. We call $\origen(\eda)$ the \textit{initial vertex} of $\eda$ and $\final(\eda)$ the \textit{terminal vertex} of $\eda$, and we say that the initial vertex and the terminal vertex of $\eda$ are the \textit{end-points} of $\eda$.\\

A \textit{path} in $\grf$ of \textit{length} $n\in\nat$ is a sequence $\pa=\vea_0,\eda_1,\vea_1,\dots,\eda_n,\vea_n$ for $n\in\nat$ satisfying that $\eda_1,\dots,\eda_n\in\grfe$, that $\vea_0,\dots,\vea_n\in\grfv$, and that $\final(\eda_i)=\vea_i$ and $\origen(\eda_i)=\vea_{i-1}$ for all $i=1,\dots,n$, and in this case we say that $\pa$ starts on $\vea_0$ and ends on $\vea_n$. We say that such a path $\pa$ is \textit{reduced} if $\eda_i\neq(\eda_{i+1})\inv$ for all $i=1,\dots, n-1$. Given a path of the form $\pa=\vea_0,\eda_1,\vea_1,\dots,\eda_n,\vea_n$, we can define its \textit{inverse path} as the path $\pa\inv=\vea_n,\eda_n\inv,\vea_{n-1},\dots,\eda_1\inv,\vea_0$. Given two paths of the form $\pa_1=\vea_0,\eda_1,\vea_1,\dots,\eda_n,\vea_n$ and $\pa_2=\vea_n,\eda_{n+1},\vea_{n+1},\dots,\eda_{n+m},\vea_{n+m}$ (the second ends on the same vertex where the first starts), we can define its \textit{concatenation} as the path $\pa_1\pa_2=\vea_0,\eda_1,\vea_1,\dots,\eda_n,\vea_n,\eda_{n+1},\vea_{n+1},\dots,\eda_{n+m},\vea_{n+m}$.\\

Let $\grf$ be a graph. A \textit{subgraph} of $\grf$ is a graph such that its set of vertices is a subset of $\grfv$, its set of edges is a subset of $\grfe$, its inverse-edge map is the restriction to its set of vertices of the inverse-edge map of $\grf$, its initial-vertex map is the restriction to its set of vertices of the initial-vertex map of $\grf$ and its terminal-vertex map is the restriction to its set of vertices of the terminal-vertex map of $\grf$. We say that $\grf$ is \textit{connected} if for every $\vea,\veb\in\grfv$ there exists a path that starts on $\vea$ and ends on $\veb$. A \textit{tree} is a graph $\arbol$ satisfying that for every $\vea,\veb\in\vertices\arbol$ there exists a unique reduced path that starts on $\vea$ and ends on $\veb$, and such a path is called the $\arbol$\textit{-geodesic} from $\vea$ to $\veb$, and it is denoted by $[\vea,\veb]_\arbol$. In particular, every tree is connected as a graph.\\

Let $\gr$ be a group acting on a set $A$. Given $a\in A$, we define the \textit{orbit} of $a$ under $\gr$ as the set $\gr a=\{\ga a \tq \ga\in\gr\}$. The set of orbits determines a partition on $A$. Given $a\in A$, we denote by $\stab_\gr(a)=\{\ga\in\gr\tq \ga a = a\}$ the \textit{stabilizer} of $a$. We say that a group $\gr$ \textit{acts on the graph} $\grf$ if $\gr$ acts on the sets $\grfv$ and $\grfe$ in such a way the maps $\df{\inv}{\grfe}{\grfe}$, $\df{\origen}{\grfe}{\grfv}$ and $\df{\final}{\grfe}{\grfv}$ are $\gr$-invariant for the action, that is, we have $(\ga\eda)\inv=\ga(\eda\inv)$, $\origen(\ga\eda)=\ga\left(\origen(\eda)\right)$ and $\final(\ga\eda)=\ga\left(\final(\eda)\right)$ for all $\ga\in\gr$ and all $\eda\in\grfe$. When we say that $\gr$ acts on a graph $\grf$, we will also assume that \textit{this action has no inversions}, meaning that there are not $\ga\in\gr$ and $\eda\in \grfe$ such that $\ga\eda=\eda\inv$.\\

Let $\gr$ be a group acting on a graph $\grf$. We define $\gr \backslash \grf$ to be the graph whose set of vertices and edges is given by $$\vertices(\gr \backslash \grf)=\{\gr \vea \tq \vea\in\grfv \} \qquad \edges(\gr \backslash \grf)=\{\gr \eda \tq \eda\in\grfe \},$$ whose inverse-edge map $\df{\inv}{\edges(\gr \backslash \grf)}{\edges(\gr \backslash \grf)}$ is defined by $(\gr \eda)\inv= \gr \left(\eda\inv\right)$, whose initial-vertex map $\df{\origen}{\edges(\gr \backslash \grf)}{\vertices(\gr \backslash \grf)}$ is defined by $\origen(\gr \eda)= \gr \left(\origen(\eda)\right)$ and whose terminal-vertex map $\df{\final}{\edges(\gr \backslash \grf)}{\vertices(\gr \backslash \grf)}$ is defined by $\final(\gr \eda)= \gr \left(\final(\eda)\right)$. Given $\grfb$ a subgraph of $\grf$, we define the \textit{projection} of $\grfb$ on $\gr \backslash \grf$ to be the subgraph of $\gr \backslash \grf$ whose set of vertices is  $\{\gr \vea \tq \vea\in\grfbv \}$ and whose set of edges is $\{\gr \eda \tq \eda\in\grfbe \}$. Moreover, for $\vea\in\grfbv$, we say that $\gr\vea$ is the projection of $\vea$ on $\gr \backslash \grf$, and for $\eda\in\grfbv$, we say that $\gr\eda$ is the projection of $\eda$ on $\gr \backslash \grf$.\\

We now use some notions from \cite{serre} that are standard in Bass-Serre theory. We are not writing down the concrete definitions as we will not need them. Instead of that, we provide the needed conclusions with the corresponding references. For the definition of graph of groups see Definition 8 from \cite{serre}. For different equivalent definitions of the fundamental group of a graph of groups see Section 5.1 from \cite{serre}. For the definition of the universal covering relative to a graph of groups see Section 5.3 from \cite{serre}. Notice that, on that reference, on the definition of the universal covering relative to a graph of groups, there is an implicit choice of a maximal tree of the underlying graph of the graph of groups, but can skip that since we are going to consider a tree as such underlying graph, so the only maximal tree is itself.

\begin{nota}
Through all the section, $\gr$ and $\grb$ will denote two non-trivial groups, and $\ce$ will denote a subgroup of $\gr*\grb$.
\end{nota}

Consider the graph of groups $(\Phi, \grf)$ whose underlying graph is $$\grf= \veagr\stackrel{\edagr}{\longrightarrow} \veacero\stackrel{\edagrb}{\longleftarrow} \veagrb$$ (we don't write the corresponding $\edagr\inv$ and $\edagrb\inv$) and that is given by $\Phi_{\veagr}=\gr$, $\Phi_{\veacero}=\{1\}$, $\Phi_{\veagrb}=\grb$, $\Phi_{\edagr}=\{1\}$ and $\Phi_{\edagrb}=\{1\}$. In Section 5.1 of \cite{serre}, at the definition b) of the fundamental group, it is provided a definition of the fundamental group using generations and relations. From this presentation, it is direct that the fundamental group of the graph of groups $(\Phi, \grf)$ is $\gr*\grb$. Consider $\arbol$ the universal covering relative to the graph of groups $(\Phi, \grf)$. By looking at the definition of the graph given in Section 5.3 of \cite{serre} and by changing the notation in order to be more clear, we can describe $\arbol$ in terms of the way $\gr*\grb$ acts as follows: the set of edges and the set of vertices are given by 
\begin{equation*}
\begin{split}
&\vertices \arbol=(\gr*\grb)\veagr\sqcup(\gr*\grb)\veacero\sqcup(\gr*\grb)\veagrb \quad \mbox{ and }\quad \edges\arbol=\edges\arbol^+\sqcup\edges\arbol^-\\
& \mbox{ where }\;\edges\arbol^+=(\gr*\grb)\edagr\sqcup(\gr*\grb)\edagrb \;\mbox{ and }\; \edges\arbol^-=(\gr*\grb)\edagr\inv\sqcup(\gr*\grb)\edagrb\inv
\end{split}
\end{equation*}
and the structure of the graph is determined by
\begin{equation*}
\begin{split}
&\gx_1\veagr=\gx_2\veagr \;\Longleftrightarrow\; \gx_1\gr=\gx_2\gr;\quad \gx_1\veacero=\gx_2\veacero \;\Longleftrightarrow\; \gx_1=\gx_2;\quad \gx_1\veagrb=\gx_2\veagrb \;\Longleftrightarrow\; \gx_1\grb=\gx_2\grb;\\
&\gx_1\edagr=\gx_2\edagr \;\Longleftrightarrow\; \gx_1=\gx_2;\quad \gx_1\edagrb=\gx_2\edagrb \;\Longleftrightarrow\; \gx_1=\gx_2;\\
&(\gx\edagr)\inv=\gx\edagr\inv;\quad (\gx\edagrb)\inv=\gx\edagrb\inv;\\
&\origen(\gx\edagr)=\gx\veagr;\quad \final(\gx\edagr)=\gx\veacero;\quad \origen(\gx\edagrb)=\gx\veagrb;\quad  \final(\gx\edagrb)=\gx\veacero.
\end{split}
\end{equation*}
The graph $\arbol$ is a tree as a consequence of Theorem 12 from \cite{serre}.\\

We can identify elements of $\gr*\grb$ with paths in $\arbol$. Notice that the action of $\gr*\grb$ on the orbit of $1\veacero$ is free. In particular, given $\gx\in\gr*\grb$, we can identify $\gx$ with a path in $\arbol$ by considering $[1\veacero,\gx\veacero]_\arbol$, which is the unique reduced path in $\arbol$ starting on $1\veacero$ and ending on $\gx\veacero$. We now describe $[1\veacero,\gx\veacero]_\arbol$. Write $\gx=\gx_1\dots\gx_n$ as a reduced word in $\gr*\grb$ (so that each $\gx_i$ is a non trivial element of $\gr$ or $\grb$ in such a way they alternate when increasing $i$ by one). The path $[1\veacero,\gx\veacero]_\arbol$ is of the form $$[1\veacero,\gx\veacero]_\arbol=\vea_0,\eda_1,\vea_1,\dots,\eda_{2n},\vea_{2n}.$$ We have $\vea_0=1\veacero$ and $\vea_{2i}=(\gx_1\dots\gx_i)\veacero$ for all $i=1,\dots, n$. Also, for each $i=1,\dots, n$, either $\eda_{2i-1}=(\gx_1\dots\gx_{i-1})\edagr\inv$, $\vea_{2i-1}=(\gx_1\dots\gx_{i-1})\veagr$ and $\eda_{2i-1}=(\gx_1\dots\gx_i)\edagr$ (this happens when $\gx_i\in\gr$), or $\eda_{2i-1}=(\gx_1\dots\gx_{i-1})\edagrb\inv$, $\vea_{2i-1}=(\gx_1\dots\gx_{i-1})\veagrb$ and $\eda_{2i-1}=(\gx_1\dots\gx_i)\edagrb$ (this happens when $\gx_i\in\grb$).

\begin{obs}\label{obs bs camino factores}
Consider $\gx\in\gr*\grb$. From the description of $[1\veacero,\gx\veacero]_\arbol$ we deduce that the following fact holds: given $\gr_0$ subgroup of $\gr$ and given $\grb_0$ subgroup of $\grb$ then 
\begin{equation*}
\begin{split}
&\gx\in \gen{\gr_0,\grb_0}\;\mbox{ and }\;[1\veacero,\gx\veacero]_\arbol=\vea_0,\eda_1,\vea_1,\dots,\eda_m,\vea_m \quad\mbox{ implies } \\
& \eda_i\in \gen{\gr_0,\grb_0}\edagr\cup\gen{\gr_0,\grb_0}\edagr\inv\cup\gen{\gr_0,\grb_0}\edagrb\cup\gen{\gr_0,\grb_0}\edagrb\inv \;\mbox{ for all } i=1,\dots,m.
\end{split}
\end{equation*}
\end{obs}

Since $\gr*\grb$ acts on $\arbol$, we can consider the induced action of $\ce$ on $\arbol$. Then, we can consider the graph $\overline{\arbol}_\ce=\ce \backslash \arbol$. Its set of vertices and edges are described by 
\begin{equation*}
\begin{split}
&\vertices \overline{\arbol}_\ce=\ce \backslash(\gr*\grb)\veagr\sqcup\ce \backslash(\gr*\grb)\veacero\sqcup\ce \backslash(\gr*\grb)\veagrb \quad \mbox{ and } \\
&\edges\overline{\arbol}_\ce=\ce \backslash(\gr*\grb)\edagr\sqcup\ce \backslash(\gr*\grb)\edagr\inv\sqcup\ce \backslash(\gr*\grb)\edagrb\sqcup\ce \backslash(\gr*\grb)\edagrb\inv,
\end{split}
\end{equation*}
where the inverse edge and the two endpoints of a given edge are the induced by the ones from $\arbol$.

\begin{obs}\label{obs bs estabilizadores descripcion}
Take $\gx\in\gr*\grb$. Notice that $$\stab_\ce(\gx\veagr)=\{\gy\in\ce\tq \gy\gx\veagr=\gx\veagr \}=\{\gy\in\ce\tq \gx\inv\gy\gx\veagr=\veagr \}=\{\gy\in\ce\tq \gx\inv\gy\gx\in\gr \}=\ce\cap(\gx\gr\gx\inv).$$ We can do analogous computations for the other vertices and edges so we obtain
\begin{equation*}
\begin{split}
&\stab_\ce(\gx\veagr)=\ce\cap(\gx\gr\gx\inv) \; \mbox{ , } \stab_\ce(\gx\veacero)=\{1\} \; \mbox{ and }\; \stab_\ce(\gx\veagrb)=\ce\cap(\gx\grb\gx\inv)\; \mbox{ on vertices, and }\\
& \stab_\ce(\gx\edagr)=\{1\} \; \mbox{ and }  \stab_\ce(\gx\edagrb)=\{1\} \mbox{ on edges.}
\end{split}
\end{equation*}
\end{obs}

\begin{obs}\label{obs bs estabilizadores triviales}
Take $\gx\in\gr*\grb$. We have that
\begin{equation*}
\stab_\ce(\gx\veagr)=\{1\}\quad \mbox{ if and only if } \quad \ce\gx\ga\neq\ce\gx \;\mbox{ for all }\ga\in\gr\settminus\{1\}
\end{equation*}
and 
\begin{equation*}
\stab_\ce(\gx\veagrb)=\{1\}\quad \mbox{ if and only if } \quad \ce\gx\gb\neq\ce\gx \;\mbox{ for all }\gb\in\grb\settminus\{1\}.
\end{equation*}
It follwos from $\stab_\ce(\gx\veagr)=\ce\cap(\gx\gr\gx\inv)$ and $\stab_\ce(\gx\veagr)=\ce\cap(\gx\gr\gx\inv)$ (done in the previous remark).
\end{obs}

Let $\cete(\ce)$ be the sum of the number of vertices of $\overline{\arbol}_\ce$ of the form $\ce\gx\veagr$ such that $\stab_\ce(\gx\veagr)\neq\{1\}$ plus the number of vertices of $\overline{\arbol}_\ce$ of the form $\ce\gx\veagrb$ such that $\stab_\ce(\gx\veagrb)\neq\{1\}$. Notice that this number is well defined as it does not depend on the chosen representatives by \Cref{obs bs estabilizadores triviales}. Let $\rank(\overline{\arbol}_\ce)$ be the number of edges of $\overline{\arbol}_\ce$ outside of a maximal subtree (this is the same as the free rank of the fundamental group of the graph $\overline{\arbol}_\ce$, which it is known that is a free group). Then, the Kurosh rank of $\ce$ with respect to $\gr*\grb$ is equal to $\cete(\ce)+\rank(\overline{\arbol}_\ce)$. This is a consequence of the proof of Kurosh's theorem, see Theorem 14 of \cite{serre}.

\begin{nota}\label{notacion bs core}
Define $\vebase=\ce 1\veacero$. Define the \textit{fundamental core} of $\overline{\arbol}_\ce$ to be the smallest connected subgraph of $\overline{\arbol}_\ce$ containing all reduced paths on $\overline{\arbol}_\ce$ simultaneously starting and ending at $\vebase$ and containing all vertices of $\overline{\arbol}_\ce$ of the form $\ce\gx\veagr$ such that $\stab_\ce(\gx\veagr)\neq\{1\}$ and all vertices of $\overline{\arbol}_\ce$ of the form $\ce\gx\veagrb$ such that $\stab_\ce(\gx\veagrb)\neq\{1\}$.
\end{nota}

\begin{obs}\label{obs bs rango de kurosh finito implica grafo finito}
Let $\ce$ be a subgroup of $\gr*\grb$. If $\ce$ has finite Kurosh rank with respect to $\gr*\grb$, then the fundamental core of $\overline{\arbol}_\ce$ is a finite graph. This is direct from the definitions. 
\end{obs}

\begin{nota}\label{notacion bs camino en cociente}
Let $\ce$ be a subgroup of $\gr*\grb$. For each $\gx\in\gr*\grb$ we define the path $\camaso(\gx)$ in $\overline{\arbol}_\ce$ as the projection of $[1\veacero,\gx\veacero]_\arbol$ on $\overline{\arbol}_\ce$, that is: if $[1\veacero,\gx\veacero]_\arbol=\vea_0,\eda_1,\vea_1,\dots,\eda_{2n},\vea_{2n}$ then $\camaso(\gx)=\ce\vea_0,\ce\eda_1,\ce\vea_1,\dots,\ce\eda_{2n},\ce\vea_{2n}$. Notice that $\camaso(\gx)$ is a path having $\vebase=\ce 1\veacero$ as initial vertex and having $\ce \gx\veacero$ as terminal vertex. Notice that $\camaso(\gx)$ always starts on $\vebase$, and that $\camaso(\gx)$ ends on $\vebase$ if and only if $\gx\in\ce$.
\end{nota}

\begin{lema}\label{lema bs camino en el core}
Let $\ce$ be a subgroup of $\gr*\grb$. If $\gx\in\ce$, then $\camaso(\gx)$ is a path on the fundamental core of $\overline{\arbol}_\ce$. In particular, if $\gx\in\ce$, then every edge and every vertex appearing on the path $[1\veacero,\gx\veacero]_\arbol$, when projected on $\overline{\arbol}_\ce$, give us an edge or a vertex that belongs to the fundamental core of $\overline{\arbol}_\ce$.
\begin{proof}
Assume that $\gx\in\ce$. In particular, $\camaso(\gx)$ is a closed path simultaneously starting and ending at $\vebase$. If $\camaso(\gx)$ is reduced, then it is on the fundamental core of $\overline{\arbol}_\ce$ by definition. If it is not reduced, we can take the minimum $i$ such that $\ce\eda_i,\ce\vea_i,\ce\eda_i\inv$, and this condition implies that $\stab_\ce(\vea_i)\neq\{1\}$. Then, by definition, the subpath of $\camaso(\gx)$ from the startng point until $\ce\vea_i$. By repeating the argument with next vertices where $\camaso(\gx)$ is not reduced, we prove the lemma. 
\end{proof}
\end{lema}

We start by proving Fact II. This will be a corollary of next result, which is stated considering fixed generating sets of $\gr$ and $\grb$ because we need it in order to derive the wanted fact.

\begin{prop}\label{proposicion general lambda g y lambda h}
Let $\gr$ and $\grb$ two non-trivial groups. Let $\ce$ be a subgroup of $\gr*\grb$ having finite Kurosh rank with respect to $\gr*\grb$. Let $\equisg$ be a generating set of $\gr$ and $\equish$ be a generating set of $\grb$. Define $\Xi\subseteq \ce\backslash (\gr*\grb)$ to be the set of cosets $\ce \lambda$ such that there exists $\ga\in \equisg$ and $\gb\in \equish$ satisfying $\ce\lambda \ga\neq \ce\lambda$, $\ce \lambda \gb\neq \ce\lambda$ and $\ce \lambda \ga^{\varepsilon_\gr}\neq \ce \lambda \gb^{\varepsilon_\grb}$ for every $\varepsilon_\gr,\varepsilon_\grb\in\{\pm 1\}$. Then, $\left(\ce\backslash (\gr*\grb)\right)\settminus \Xi$ is a finite set.
\begin{proof}
By \Cref{obs bs rango de kurosh finito implica grafo finito}, since $\ce$ has finite Kurosh rank with respect to $\gr*\grb$, we deduce that the fundamental core of $\overline{\arbol}_\ce$ is finite. If we show that all vertices of the form $\ce\lambda\veacero$ for $\ce\lambda\in\left(\ce\backslash (\gr*\grb)\right)\settminus \Xi$ are vertices of the fundamental core of $\overline{\arbol}_\ce$, we can conclude that $\left(\ce\backslash (\gr*\grb)\right)\settminus \Xi$ is finite.\\

Indeed, assume that $\ce\lambda\veacero$ for $\ce\lambda\in\ce\backslash (\gr*\grb)$ is not on the fundamental core of $\overline{\arbol}_\ce$. We want to prove that $\ce\lambda$ belongs to $\Xi$. Consider the edges $\ce\lambda\edagr$ and $\ce\lambda\edagrb$ of $\overline{\arbol}_\ce$ and the vertices $\ce\lambda\veagr$ and $\ce\lambda\veagrb$ of $\overline{\arbol}_\ce$ that satisfies $\origen(\ce\lambda\edagr)=\ce\lambda\veagr$, $\origen(\ce\lambda\edagrb)=\ce\lambda\veagrb$ and $\final(\ce\lambda\edagr)=\final(\ce\lambda\edagrb)=\ce\lambda\veacero$. In other words, we have $$\ce\lambda\veagr\stackrel{\ce\lambda\edagr}{\longrightarrow} \ce\lambda\veacero\stackrel{\ce\lambda\edagrb}{\longleftarrow} \ce\lambda\veagrb.$$ We remark that either $\ce\lambda\veagr$ or $\ce\lambda\veagrb$ is not on the fundamental core of $\overline{\arbol}_\ce$. Indeed, if both $\ce\lambda\veagr$ and $\ce\lambda\veagrb$ are on the fundamental core of $\overline{\arbol}_\ce$, since the path $\ce\lambda\veagr,\ce\lambda\edagrb,\ce\lambda\veacero,(\ce\lambda\edagrb)\inv,\ce\lambda\veagrb$ is a reduced path starting on $\ce\lambda\veagr$ and ending on $\ce\lambda\veagrb$, it would imply that $\ce\lambda\veacero$ is on the fundamental core of $\overline{\arbol}_\ce$, a contradiction with the fact that it is not on the fundamental core of $\overline{\arbol}_\ce$. Therefore, we know that either $\ce\lambda\veagr$ or $\ce\lambda\veagrb$ is not on the fundamental core of $\overline{\arbol}_\ce$. We consider the case where $\ce\lambda\veagr$ is not on the fundamental core of $\overline{\arbol}_\ce$, as the case where $\ce\lambda\veagrb$ is not on the fundamental core of $\overline{\arbol}_\ce$ is analogous.\\

Assume that $\ce\lambda\veagr$ is not on the fundamental core of $\overline{\arbol}_\ce$. Fix a nontrivial $\ga\in \equisg$. Since $\ce\lambda\veagr$ is not on the fundamental core of $\overline{\arbol}_\ce$, we have $\stab_\ce(\lambda\veagrb)=\{1\}$. Then, by \Cref{obs bs estabilizadores triviales}, we deduce that 
\begin{equation}\label{proposicion general lambda g y lambda h eq conclusion previa 1}
\ce\lambda\ga\neq\ce\lambda \;\mbox{ for all }\; \ga\in \equisg\settminus\{1\}.
\end{equation}
We now prove that 
\begin{equation}\label{proposicion general lambda g y lambda h eq conclusion previa 2}
\mbox{There exists }\;\gb\in \equish\;\mbox{ such that }\; \ce\lambda\gb\neq\ce\lambda.
\end{equation}
We do it by contradiction, so assume that $\ce\lambda\gb=\ce\lambda$ for all $\gb\in \equish$. Since $\equish$ is a generating set of $\grb$, we deduce that
\begin{equation}\label{proposicion general lambda g y lambda h eq auxiliar a conclusion previa 2}
\ce\lambda\gb=\ce\lambda \;\mbox{ for all }\; \gb\in \grb.
\end{equation}
This implies, by the structure of $\overline{\arbol}_\ce$, that $\ce\lambda\edagrb$ is the unique edge of $\overline{\arbol}_\ce$ having $\ce\lambda\veagrb$ as initial vertex. \Cref{proposicion general lambda g y lambda h eq auxiliar a conclusion previa 2} implies that $\ce\lambda\veagrb$ is on the fundamental core of $\overline{\arbol}_\ce$. In particular, any path starting on $\ce\lambda\veagrb$ and ending on $\vebase$ pass through the vertex $\ce\lambda\veacero$. Hence, $\ce\lambda\veacero$ is on the fundamental core of $\overline{\arbol}_\ce$, a contradiction. Therefore, we have shown \Cref{proposicion general lambda g y lambda h eq conclusion previa 2}.\\

Fix $\ga\in \equisg\settminus\{1\}$ satisfying that
$\ce\lambda\ga^{\varepsilon_\gr}\neq\ce\lambda$ for all $\varepsilon_\gr\in\{\pm 1\}$, which is possible by \Cref{proposicion general lambda g y lambda h eq conclusion previa 1} since $\gr$ is non-trivial. Fix $\gb\in \equish$ satisfying $\ce\lambda\gb^{\varepsilon_\grb}\neq\ce\lambda$ for all $\varepsilon_\grb\in\{\pm 1\}$, which is possible by \Cref{proposicion general lambda g y lambda h eq conclusion previa 2}. Fix $\varepsilon_\gr,\varepsilon_\grb\in\{\pm 1\}$. In particular, we have that $\ce\lambda \ga^{\varepsilon_\gr}\neq \ce\lambda$ and $\ce \lambda \gb^{\varepsilon_\grb}\neq \ce\lambda$. We can finish the proof if we show that 
\begin{equation}\label{proposicion general lambda g y lambda h eq conclusion previa 3}
\ce \lambda \ga^{\varepsilon_\gr}\neq \ce \lambda \gb^{\varepsilon_\grb}.
\end{equation}
We do it by contradiction, so assume that $\ce\lambda\ga^{\varepsilon_\gr}=\ce\lambda\gb^{\varepsilon_\grb}$. In particular, $\ce\lambda\ga^{\varepsilon_\gr}\veacero=\ce\lambda\gb^{\varepsilon_\grb}\veacero$. Then, we can consider the path $$\pa=\ce\lambda\ga^{\varepsilon_\gr}\veacero, \ce\lambda\ga^{\varepsilon_\gr}\edagr\inv, \ce\lambda\veagr, \ce\lambda\edagr, \ce\lambda\veacero, \ce\lambda\edagrb\inv, \ce\lambda\veagrb, \ce\lambda\gb^{\varepsilon_\grb}\edagrb, \ce\lambda\gb^{\varepsilon_\grb}\veacero,$$ which is a path starting and ending on $\ce\lambda\ga^{\varepsilon_\gr}\veacero=\ce\lambda\gb^{\varepsilon_\grb}\veacero$, and notice that it is reduced since we know that $\ce\lambda\ga^{\varepsilon_\gr}\neq\ce\lambda$ and $\ce\lambda\gb^{\varepsilon_\grb}\neq\ce\lambda$. In particular, and since $\ce\lambda\veacero$ is a vertex appearing on $\pa$, we can use this fact to build a reduced path having $\ce\lambda\veacero$ as a vertex and starting and ending on $\vebase$. Hence, by definition of the fundamental core of $\overline{\arbol}_\ce$, we deduce that $\ce\lambda\veacero$ is on the fundamental core of $\overline{\arbol}_\ce$. This contradicts the fact that $\ce\lambda\veagr$ is not on the fundamental core of $\overline{\arbol}_\ce$ by hypothesis. Therefore, we have shown \Cref{proposicion general lambda g y lambda h eq conclusion previa 3}, as we wanted.
\end{proof}
\end{prop}

As a corollary, we obtain Fact II by taking inverses on the prevous result, so that right cosets are transformed into left cosets.

\begin{coro}\label{corolario final lambda g y lambda h}
Let $\gr$ and $\grb$ two non-trivial groups. Let $\ce$ be a subgroup of $\gr*\grb$ having finite Kurosh rank with respect to $\gr*\grb$.  Let $\equisg$ be a generating set of $\gr$ and $\equish$ be a generating set of $\grb$. Let $(\lambda_l)_{n\in\nat}$ be a sequence of elements in $\gr*\grb$ such that $\lambda_l \ce\neq\lambda_m \ce$ for all $l\neq m$. Then, there exists $n\in\nat$ such that for all $k\geq n$ there are $\ga\in \equisg$ and $\gb\in \equish$ satisfying $\ga \lambda_k \ce\neq\lambda_k \ce$, $\gb \lambda_k \ce\neq\lambda_k \ce$ and $\ga^{\varepsilon_\gr} \lambda_k \ce\neq \gb^{\varepsilon_\grb}\lambda_k \ce$ for all $\varepsilon_\gr,\varepsilon_\grb\in\{\pm 1\}$.
\end{coro}

We conclude with the main result of the section, which corresponds to Fact I. As we said at the beginning of the section, as far as the author knows, this conclusion has not been considered in the literature. In \cite{kurosh-antolin} it is proven that the reduced Kurosh rank of the intersection of two subgroups of a free product of left-orderable groups is bounded above by the product of the reduced Kurosh ranks of each subgroup. We provide the same conclusion for  a certain case where one of the intersecting subgroups has Kurosh rank $2$ but dropping left-orderable requirement. Notice that this property is a generalization of strengthened Hanna Neumann Conjecture.  Related to this, it is remarkable that Tardos proved strengthened Hanna Neumann Conjecture on free groups for the case where one of the subgroups has free rank $2$ (\cite{tardos}).

\begin{teo}\label{bass serre para interseccion}
Let $\gr$ and $\grb$ two non-trivial groups. Let $\ce$ be a proper subgroup of $\gr*\grb$ such that $\ce$ has finite Kurosh rank with respect to $\gr*\grb$. Then, there exists finite sets $\fino,\{\zeta_\indi\}_{\indi\in\fini},\{\alpha_\indj\}_{\indj\in\finj},\{\beta_\indk\}_{\indk\in\fink}\subseteq \gr*\grb$ satisfying the following property:\\
For every subgroup $\gr_0$ of $\gr$ and every subgroup $\grb_0$ of $\grb$ satisfying $\fino\subseteq \gen{\gr_0,\grb_0}$ there exists $\gr_\indj$ subgroup of $\gr_0$ for each $\indj\in\finj$ and there exists $\gr_\indk$ subgroup of $\grb_0$ for each $\indk\in\fink$ such that 
\begin{equation}\label{bass serre para interseccion eq enunciado}
\ce\cap\gen{\gr_0,\grb_0}=\gen{\{\zeta_\indi \tq \indi\in\fini\}\cup \left(\bigcup_{\indj\in\finj}\alpha_\indj\gr_\indj\alpha_\indj\inv \right)\cup \left(\bigcup_{\indk\in\fink}\beta_\indk\grb_\indk\beta_\indk\inv \right)}.
\end{equation}

\begin{proof}

Since $\ce$ has finite Kurosh rank with respect to $\gr*\grb$, applying \Cref{obs bs rango de kurosh finito implica grafo finito}, we deduce that the fundamental core of $\overline{\arbol}_\ce$ is a finite subgraph of $\overline{\arbol}_\ce$. Then, we can take $\overline{\arbmax}$ a maximal tree of the fundamental core of $\overline{\arbol}_\ce$, which is finite and satisfies that every vertex of the fundamental core of $\overline{\arbol}_\ce$ is a vertex of $\overline{\arbmax}$. Since $\overline{\arbmax}$ is finite, we can recursively build $\arbmax$ a subgraph of $\arbol$ that is a finite tree containing $1\veacero$ in such a way the projection of $\arbmax$ on $\overline{\arbol}_\ce$ is $\overline{\arbmax}$ and every vertex and every edge of $\overline{\arbmax}$ is the projection of a unique vertex or edge of $\arbmax$. Notice that every vertex of the fundamental core of $\overline{\arbol}_\ce$ is the projection of a unique vertex of $\arbmax$. We now construct a fundamental tranvsersal $\aristasfuera$ (see \cite{actingongraphs-dicksdunwoody} for a definition). We give all the details to make the exposition sefl-contained. Since $\overline{\arbmax}$ is finite, we can choose $\aristasfuera$ a finite subset of $\edges\arbol^+\settminus\edges\arbmax$ such that every $\edb\in\aristasfuera$ satisfies that $\origen(\edb)\in\vertices\arbmax$ and such that every edge of the fundamental core of $\overline{\arbol}_\ce$ not in $\overline{\arbmax}$ is the projection of a unique element in $\aristasfuera\sqcup\aristasfuera\inv$. Define $\aristasfuera^{\pm}=\aristasfuera\sqcup\aristasfuera\inv$. Hence, every vertex of the fundamental core of $\overline{\arbol}_\ce$ is the projection of a unique element of $\vertices\arbol$ and every edge of the fundamental core of $\overline{\arbol}_\ce$ is the projection of a unique element of $\edges\arbol\sqcup\aristasfuera^{\pm}$.\\

We now provide some definitions:
\begin{itemize}
\item Define $\fini=\aristasfuera$, which is a finite set.
\item Define $\finj$ to be the set of all $\gx\veagr$ for $\gx\in\gr*\grb$ satisfying that $\gx\veagr\in\vertices\arbmax$ and that $\stab_\ce(\gx\veagr)\neq\{1\}$. The fact that $\arbmax$ is finite implies that $\finj$ is finite.
\item Define $\fink$ to be the set of all $\gx\veagrb$ for $\gx\in\gr*\grb$ satisfying that $\gx\veagrb\in\vertices\arbmax$ and that $\stab_\ce(\gx\veagr)\neq\{1\}$. The fact that $\arbmax$ is finite implies that $\fink$ is finite.
\item For each $\indi\in\fini$, there exists an element of $\vertices\arbmax$ whose projection is $\ce\final(\indi)$, and this implies that $\zeta_\indi\final(\indi)\in\vertices\arbmax$ for some $\zeta_\indi\in\ce$. Hence, for each $\indi\in\fini$, we define $\zeta_\indi\in\ce$ such that $\zeta_\indi\final(\indi)\in\vertices\arbmax$.
\item For each $\indj\in\finj$ define $\alpha_\indj\in\gr*\grb$ to be such that $\origen(\alpha_\indj\edagr)=\indj$ and such that $\alpha_\indj\edagr$ is an edge appearing on the the $\arbol$-geodesic from $1\veacero$ to $\indj$.
\item For each $\indk\in\fink$ define $\beta_\indk\in\gr*\grb$ to be such that $\origen(\beta_\indk\edagrb)=\indk$ and such that $\beta_\indk\edagrb$ is an edge appearing on the the $\arbol$-geodesic from $1\veacero$ to $\indk$.
\end{itemize}

Note that, in the previous definitions, the elements $\zeta_\indi$, $\alpha_\indj$ and $\beta_\indk$ are the unique elements having the corresponding properties. This is a consequence of the fact that the action of $\gr*\grb$ on $\arbol$ is free on edges and on vertices in the orbit of $\veacero$.\\

Define the finite set $$\fino=\{\zeta_\indi \tq \indi\in\fini\}\cup \{\alpha_\indj \tq \indj\in\finj\}\cup \{\beta_\indk \tq \indk\in\fink\}\cup \{\gx\in\gr*\grb\tq \gx\edagr\in \edges\arbmax\cup\aristasfuera^{\pm}\;\mbox{ or }\;\gx\edagrb\in \edges\arbmax\cup\aristasfuera^{\pm}\}$$ We are going to prove that the theorem holds for these subsets we have defined. Consider $\gr_0$ a subgroup of $\gr$ and $\grb_0$ a subgroup of $\grb$ satisfying that
\begin{equation*}
\fino\subseteq \gen{\gr_0,\grb_0}.
\end{equation*}
Define
\begin{equation*}
\gr_\indj=\left(\alpha_\indj\inv\stab_\ce(\alpha_\indj\edagr)\alpha_\indj\right)\cap\gr_0 \;\mbox{ for }\, \indj\in\finj \quad\mbox{ and }\quad \grb_\indk=\left(\beta_\indk\inv\stab_\ce(\beta_\indk\edagrb)\beta_\indk\right)\cap\grb_0  \;\mbox{ for }\, \indk\in\fink.
\end{equation*}
This, by \Cref{obs bs estabilizadores descripcion}, is equivalent to
\begin{equation}\label{bass serre para interseccion eq subgrupos de la descomposicion de kurosh}
\gr_\indj=\left(\alpha_\indj\inv\ce\alpha_\indj\right)\cap\gr_0 \;\mbox{ for }\, \indj\in\finj \quad\mbox{ and }\quad \grb_\indk=\left(\beta_\indk\inv\ce\beta_\indk\right)\cap\grb_0  \;\mbox{ for }\, \indk\in\fink.
\end{equation}
Define 
\begin{equation*}
\generadores=\{\zeta_\indi \tq \indi\in\fini\}\cup \left(\bigcup_{\indj\in\finj}\alpha_\indj\gr_\indj\alpha_\indj\inv \right)\cup \left(\bigcup_{\indk\in\fink}\beta_\indk\grb_\indk\beta_\indk\inv \right).
\end{equation*}
We want to prove is that $\ce\cap\gen{\gr_0,\grb_0}=\gen{\generadores}$. We first prove that 
\begin{equation}\label{bass serre para interseccion eq generadores contenidos en la interseccion}
\gen{\generadores}\subseteq \ce\cap\gen{\gr_0,\grb_0}.
\end{equation}
In order to do that, we check that $\generadores\subseteq \ce\cap\gen{\gr_0,\grb_0}$. The fact that $\zeta_\indi\in\ce\cap\gen{\gr_0,\grb_0}$ for all $\indi\in\fini$ follows by the fact that $\zeta_\indi\in\ce$ for all $\indi\in\fini$ by definition and that $\zeta_\indi\in\fino\subseteq \gen{\gr_0,\grb_0}$ for all $\indi\in\fini$. The fact that, for each $\indj\in\finj$, we have $\alpha_\indj\gr_\indj\alpha_\indj\inv\subseteq\ce\cap\gen{\gr_0,\grb_0}$, follows from $\gr_\indj=\left(\alpha_\indj\inv\ce\alpha_\indj\right)\cap\gr_0$ by \Cref{bass serre para interseccion eq subgrupos de la descomposicion de kurosh} and from $\alpha_\indj\in\fino\subseteq \gen{\gr_0,\grb_0}$. The fact that $\beta_\indk\grb_\indk\beta_\indk\inv\subseteq\ce\cap\gen{\gr_0,\grb_0}$ for each $\indk\in\fink$ is analogous. That completes the proof of \Cref{bass serre para interseccion eq generadores contenidos en la interseccion}.\\

Recall that all we wanted to prove is $\ce\cap\gen{\gr_0,\grb_0}=\gen{\generadores}$. By \Cref{bass serre para interseccion eq generadores contenidos en la interseccion}, we only need to show that $\ce\cap\gen{\gr_0,\grb_0}\subseteq \gen{\generadores}$.\\

We now introduce some auxiliary definitions. A vertex $\vea$ of $\arbol$ is said to be critical if $\vea\notin\vertices\arbmax$. An edge $\eda$ of $\arbol$ is said to be critical if $\eda\notin\edges\arbmax\cup\aristasfuera^{\pm}$. Define the critical number of a path $\pa$ to be the number of critical vertices appearing on $\pa$ plus the number of critical edges appearing on $\pa$. Notice that, for any $\vea_1,\vea_2\in\vertices\arbol$, the critical number of $[\vea_1,\vea_2]_\arbol$ is smaller or equal than the critical number of any path starting on $\vea_1$ and ending on $\vea_2$. For each $\gx\in\ce\cap\gen{\gr_0,\grb_0}$, define the critical number of $\gx$ to be the critical number of $[1\veacero,\gx\veacero]_\arbol$. We are going to prove $\ce\cap\gen{\gr_0,\grb_0}\subseteq \gen{\generadores}$ by induction on the critical number.\\

Now, we are ready to do the induction. Assume that $\gx\in\ce\cap\gen{\gr_0,\grb_0}$ has critical number $0$. In particular, $\gx\veacero$ is a vertex of $\arbmax$. As $\gx\in\ce$, we deduce that $1\veacero$ and $\gx\veacero$ both project to the same vertex of $\overline{\arbmax}$. We deduce that $1\veacero=\gx\veacero$. Hence $\gx=1$, so then $\gx\in\gen{\generadores}$.\\

Before doing the inductive step, we prove a claim that we will use few times later on the inductive proof. We claim that, for each $\gx\in\ce\cap\gen{\gr_0,\grb_0}$, if 
\begin{equation}\label{bass serre para interseccion eq forma del camino 0}
\begin{split}
&[1\veacero,\gx\veacero]_\arbol=1\veacero,\eda_1,\vea_1,\dots,\eda_i,\vea_i, \eda_{i+1},\dots,\eda_n,\gx\veacero \\
 & \mbox{ where } \vea_1,\dots,\vea_{n-1}\in\vertices\arbol,\; \eda_1,\dots,\eda_{n-1}\in\edges\arbol,\; n\in\natp,
\end{split}
\end{equation}
if $\vea_i$ (respectively $\eda_i$) is critical for some $i$, and if $\gy\in\ce$ is such that $\gy\vea_i\in \vertices\arbmax$ (resp. $\gy\eda_i\in \edges\arbmax$), then $\gy\gx$ has critical number strictly smaller than the critical number of $\gx$.\\

We prove the claim. Consider $\gx\in\ce\cap\gen{\gr_0,\grb_0}$. Assume that $[1\veacero,\gx\veacero]_\arbol$ is of the form as in \Cref{bass serre para interseccion eq forma del camino 0}, that $\vea_i$ is critical for some $i$, and that $\gy\in\ce$ is such that $\gy\vea_i\in \vertices\arbmax$ (the case for edges is analogous). Recall that $[1\veacero,\gx\veacero]_\arbol$ is a reduced path on the tree $\arbol$ where $1\veacero$ is a vertex of the tree $\arbmax$. From this fact and from the properties of $\aristasfuera$, we deduce that the path $[1\veacero,\gx\veacero]_\arbol$ is a concatenation of a sequence of non critical elements and then a sequence of critical elements. Then, by erasing on that path all the elements appearing before $\vea_i$ and left-multiplying by $\gy$, we obtain a path starting on $\gy\vea_i$, ending on $\gy\gx\veacero$ and having critical number strictly smaller than the critical number of $[1\veacero,\gx\veacero]_\arbol$ (that is, the critical number of $\gx$). By concatenating $[1\veacero,\gy\vea_i]_\arbol$ with that path, since all elements of $[1\veacero,\gy\vea_i]_\arbol$ are not critical as it is a path on $\arbmax$, we obtain a path starting on $1\veacero$, ending on $\gy\gx\veacero$ and having critical number strictly smaller than  the critical number of $\gx$. Then, $[1\veacero,\gy\gx\veacero]_\arbol$ has critical number strictly smaller than the critical number of $\gx$. This proves the claim.\\

We are now ready to study the inductive case. Consider $\gx\in\ce\cap\gen{\gr_0,\grb_0}$ having critical number $s>0$. Then, $[1\veacero,\gx\veacero]_\arbol$ has some critical vertex or some critical edge. Recall that $[1\veacero,\gx\veacero]_\arbol$ is a reduced path on the tree $\arbol$ where $1\veacero$ is a vertex of the tree $\arbmax$. From this and the properties of $\aristasfuera$, we deduce that the path $[1\veacero,\gx\veacero]_\arbol$ is either a sequence of elements of $\arbmax$ (ending with a vertex) and then critical elements, or a sequence of elements of $\arbmax$ (ending with a vertex), followed by an edge of $\aristasfuera$ and then critical elements. Knowing this, we have some distinct cases. In all the cases, we are going to use the fact that all edges and vertices appearing in $[1\veacero,\gx\veacero]_\arbol$, when projected to $\overline{\arbol}_\ce$, they are all elements of the fundamental core of $\overline{\arbol}_\ce$, since $\gx\in\ce$ and as a consequence of \Cref{lema bs camino en el core}.
\begin{itemize}

\item[\textbf{Case I}] Suppose that $[1\veacero,\gx\veacero]_\arbol$ is a sequence of elements of $\arbmax$ (ending with a vertex), followed by an edge of $\aristasfuera$ and then critical elements. In other words and more formally, we can write 
\begin{equation*}
\begin{split}
&[1\veacero,\gx\veacero]_\arbol=1\veacero,\eda_1',\vea_1',\dots,\eda_{n-1}',\vea_{n-1}',\eda_n',\vea_1,\eda_1,\dots,\vea_r,\eda_r,\gx\veacero \\
& \mbox{ where } \vea_1',\dots,\vea_{n-1}'\in\vertices\arbmax,\; \eda_1',\dots,\eda_{n-1}'\in\edges\arbmax,\;\eda_n'\in\aristasfuera,\\
&\vea_1,\dots,\vea_r\in\vertices\arbol\settminus\vertices\arbmax,\; \eda_1,\dots,\eda_r\in\edges\arbol\settminus(\edges\arbmax\cup\aristasfuera^{\pm}), \; n\in\natp,\;r\in\nat,
\end{split}
\end{equation*}
where the critical number of $[1\veacero,\gx\veacero]_\arbol$ is $s=2r+1$.\\

Since $\eda_n'\in\aristasfuera$ and $\aristasfuera=\indi$, we deduce that $\zeta_{\eda_n'}\final(\eda_n')\in\vertices\arbmax$, so that $\zeta_{\eda_n'}\vea_1\in\vertices\arbmax$. Then, by the claim, we deduce that $\zeta_{\eda_n'}\gx$ has critical number strictly smaller than the critical number of $[1\veacero,\gx\veacero]_\arbol$, which is $s$.\\

We have $\zeta_{\eda_n'}\in\fino\subseteq \gen{\gr_0,\grb_0}$. As we know that $\zeta_{\eda_n'}\in\ce$ by definition of $\zeta_{\eda_n'}$, we obtain $\zeta_{\eda_n'}\in\ce\cap\gen{\gr_0,\grb_0}$. Applying that $\gx\in\ce\cap\gen{\gr_0,\grb_0}$ by hypothesis, we obtain that $\zeta_{\eda_n'}\gx\in\ce\cap\gen{\gr_0,\grb_0}$. Then, since $\zeta_{\eda_n'}\gx$ has critical number strictly smaller than $s$, by induction hypothesis, we obtain that $\zeta_{\eda_n'}\gx\in \gen{\generadores}$. As we know that $\zeta_{\eda_n'}\in \gen{\generadores}$ from the definition of $\generadores$, we conclude that $\gx\in \gen{\generadores}$. This finishes the Case I.\\

\item[\textbf{Case II}]  Suppose that the first critical element is a critical edge, and suppose that the previous vertex is an element of the form $\gy\veacero$ for some $\gy\in\gr*\grb$. Formally, this means that we can write 
\begin{equation*}
\begin{split}
&[1\veacero,\gx\veacero]_\arbol=1\veacero,\eda_1',\vea_1',\dots,\eda_n',\vea_n',\eda_1,\vea_1,\dots,\eda_r,\gx\veacero \;\mbox{ where } \vea_1',\dots,\vea_n'\in\vertices\arbmax,\; \eda_1',\dots,\eda_n'\in\edges\arbmax \\
& \vea_1,\dots,\vea_{r-1}\in\vertices\arbol\settminus\vertices\arbmax,\; \eda_1,\dots,\eda_r\in\edges\arbol\settminus(\edges\arbmax\cup\aristasfuera^{\pm}), \; \vea_n'=\gy\veacero, \;\gy\in\gr*\grb, \; n\in\nat,\;r\in\natp,
\end{split}
\end{equation*}
where the critical number of $[1\veacero,\gx\veacero]_\arbol$ is $s=2r$.\\

Notice that $\eda_1$ projects on $\overline{\arbol}_\ce$ to an edge of $\overline{\arbol}_\ce$ that is not an edge of $\overline{\arbmax}$ or of $\aristasfuera$. Then, $\eda_1$ projects on $\overline{\arbol}_\ce$ to an edge of $\aristasfuera\inv$. We deduce that there exists an element $\hat{\eda}\in\aristasfuera$ such that both $\hat{\eda}$ and $\eda_1\inv$ projects on $\overline{\arbol}_\ce$ to the same edge. Then, by looking at their final vertices and using that $\zeta_{\hat{\eda}}\in\ce$, we deduce that $\final(\eda_1\inv)=\vea_n'$ and $\zeta_{\hat{\eda}}\final(\hat{\eda})\in\vertices\arbmax$ projects on $\overline{\arbol}_\ce$ to the same vertex, so then $\vea_n'=\zeta_{\hat{\eda}}\final(\hat{\eda})$. In particular, since the action of $\ce$ is free on the orbits of $\vea_n'$ (the same orbit as $\final(\hat{\eda})$) and $\eda_1\inv$ (the same orbit as $\hat{\eda}$), we deduce that $\eda_1\inv=\zeta_{\hat{\eda}}\hat{\eda}$. Therefore, by looking at their initial vertices, we obtain $\zeta_{\hat{\eda}}\inv\vea_1\in\vertices\arbmax$. Then, by the claim, we obtain that $\zeta_{\hat{\eda}}\inv\gx$ has critical number strictly smaller than the critical number of $[1\veacero,\gx\veacero]_\arbol$, which is $s$.\\

We have $\zeta_{\hat{\eda}}\in\fino\subseteq \gen{\gr_0,\grb_0}$. As $\zeta_{\hat{\eda}}\in\ce$ by definition, we deduce $\zeta_{\hat{\eda}}\in\ce\cap\gen{\gr_0,\grb_0}$. Applying that $\gx\in\ce\cap\gen{\gr_0,\grb_0}$ by hypothesis, we obtain that $\zeta_{\hat{\eda}}\inv\gx\in\ce\cap\gen{\gr_0,\grb_0}$. Then, as $\zeta_{\hat{\eda}}\inv\gx$ has critical number strictly smaller than $s$, by induction hypothesis, we obtain that $\zeta_{\hat{\eda}}\inv\gx\in \gen{\generadores}$. As we know that $\zeta_{\hat{\eda}}\in \gen{\generadores}$ from the definition of $\generadores$, we obtain that $\gx\in \gen{\generadores}$. This concludes the Case II.\\

\item[\textbf{Case III}]  Suppose that the first apparition of a critical element in $[1\veacero,\gx\veacero]_\arbol$ is a critical edge, and suppose that the previous vertex is an element of the form $\gy\veagr$ for some $\gy\in\gr*\grb$. That is, we can write 
\begin{equation*}
\begin{split}
&[1\veacero,\gx\veacero]_\arbol=1\veacero,\eda_1',\vea_1',\dots,\eda_n',\vea_n',\eda_1,\vea_1,\dots,\eda_r,\gx\veacero \quad \mbox{ where } \vea_1',\dots,\vea_n'\in\vertices\arbmax,\; \eda_1',\dots,\eda_n'\in\edges\arbmax,\\
& \vea_1,\dots,\vea_{r-1}\in\vertices\arbol\settminus\vertices\arbmax,\; \eda_1,\dots,\eda_r\in\edges\arbol\settminus(\edges\arbmax\cup\aristasfuera^{\pm}), \; \vea_n'=\gy\veagr, \;\gy\in\gr*\grb, \; n\in\nat,\;r\in\natp,
\end{split}
\end{equation*}
where the critical number of $[1\veacero,\gx\veacero]_\arbol$ is $s$. Since $\vea_n'=\gy\veagr$ and $\vea_n'\in\vertices\arbmax$, from the definition of $\finj$, we deduce that $\vea_n'\in\finj$. In particular, we deduce that 
\begin{equation}\label{bass serre para interseccion eq descripcion 3}
\eda_n'=\alpha_{\vea_n'}\edagr\inv \; \mbox{, } \; \vea_n'=\alpha_{\vea_n'}\veagr \; \mbox{ and } \; \eda_1=\alpha_{\vea_n'}\ga\edagr \quad \mbox{ for some } \; \ga\in\gr.
\end{equation}
By definition of $\arbmax$ and $\aristasfuera$, we deduce that there is $\hat{\eda}\in \edges\arbmax\cup\aristasfuera$ such that $\hat{\eda}$ and $\eda_1$ projects on $\overline{\arbol}_\ce$ to the same edge. Then, since the initial vertex of $\eda_1$ is $\vea_n'\in\vertices\arbmax$, we deduce that the initial vertex of $\eda_1$ and the initial vertex of $\hat{\eda}$ are both equal to $\vea_n'$. Then, from \Cref{bass serre para interseccion eq descripcion 3} and from the definition of $\arbol$, obtain
\begin{equation}\label{bass serre para interseccion eq descripcion nueva 3}
\eda_n'=\alpha_{\vea_n'}\edagr\inv \; \mbox{, } \; \vea_n'=\alpha_{\vea_n'}\veagr \; \mbox{, } \; \eda_1=\alpha_{\vea_n'}\ga\edagr\; \mbox{ and } \; \hat{\eda}=\alpha_{\vea_n'}\hat{\ga}\edagr \;\; \mbox{ for some } \; \ga,\hat{\ga}\in\gr \; \mbox{ and }\; \hat{\eda}\in \edges\arbmax\cup\aristasfuera.
\end{equation}
Notice that $\alpha_{\vea_n'}\in \gen{\gr_0,\grb_0}$ and that $\alpha_{\vea_n'}\hat{\ga}\in \gen{\gr_0,\grb_0}$ since $\hat{\eda}=\hat{\ga}\alpha_{\vea_n'}\edagr$ (these facts are consequence of these elements belonging to $\fino$, which is contained in $\gen{\gr_0,\grb_0}$). Also notice that $\alpha_{\vea_n'}\ga \in \gen{\gr_0,\grb_0}$ by applying \Cref{obs bs camino factores} (which says that every vertex and edge on $[1\veacero,\gx\veacero]_\arbol$ is in $\gen{\gr_0,\grb_0}$ if we assume that $\gx\in \gen{\gr_0,\grb_0}$). From these facts we deduce that $\hat{\ga}\ga\inv,\alpha_{\vea_n'}\hat{\ga}\ga\inv\alpha_{\vea_n'}\inv \in \gen{\gr_0,\grb_0}$. Moreover, by the properties of reduced word in $\gr*\grb$, since $\gr_0$ is a subgroup of $\gr$ and $\grb_0$ is a subgroup of $\grb$ and since $\hat{\ga}\ga\inv\in\gr$, we have obtained
\begin{equation}\label{bass serre para interseccion eq el elemento esta en factores 3}
\hat{\ga}\ga\inv\in \gr_0\quad\mbox{ and }\quad \alpha_{\vea_n'}\hat{\ga}\ga\inv\alpha_{\vea_n'}\inv \in \gen{\gr_0,\grb_0}.
\end{equation}
Recall that, by definition of $\hat{\eda}$, we have that $\hat{\eda}$ and $\eda_1$ projects on $\overline{\arbol}_\ce$ to the same edge. Then, knowing that $\hat{\eda}=\alpha_{\vea_n'}\hat{\ga}\edagr$ and $\eda_1=\alpha_{\vea_n'}\ga\edagr$ from \Cref{bass serre para interseccion eq descripcion nueva 3}, we deduce that $\ce\alpha_{\vea_n'}\hat{\ga}\edagr=\ce\alpha_{\vea_n'}\ga\edagr$ as edges of $\arbol$. By definition of the edges of $\arbol$, we obtain $\ce\alpha_{\vea_n'}\hat{\ga}=\ce\alpha_{\vea_n'}\ga$. Hence, $\hat{\ga}\ga\inv\in\alpha_{\vea_n'}\inv\ce\alpha_{\vea_n'}$. Then, since $\hat{\ga}\ga\inv\in \gr_0$ by \Cref{bass serre para interseccion eq el elemento esta en factores 3}, we obtain $\hat{\ga}\ga\inv\in \left(\alpha_{\vea_n'}\inv\ce\alpha_{\vea_n'}\right)\cap\gr_0$. In other words, by definition of $\gr_{\vea_n'}$ (see \Cref{bass serre para interseccion eq subgrupos de la descomposicion de kurosh}), we have shown that $\hat{\ga}\ga\inv\in\gr_{\vea_n'}$. In particular, we obtain
\begin{equation}\label{bass serre para interseccion eq el elemento esta en gen de S 3}
\alpha_{\vea_n'}\hat{\ga}\ga\inv\alpha_{\vea_n'}\inv \in\gen{S},
\end{equation}
when applying the definition of $\generadores$.\\

Recall that $\eda_1=\alpha_{\vea_n'}\ga\edagr$ and $\hat{\eda}=\alpha_{\vea_n'}\hat{\ga}\edagr$ by \Cref{bass serre para interseccion eq descripcion nueva 3}, so we deduce that
\begin{equation*}
\left(\alpha_{\vea_n'}\hat{\ga}\ga\inv\alpha_{\vea_n'}\inv\right)\eda_1=\hat{\eda},
\end{equation*}
where $\hat{\eda}\in\vertices\arbmax$. Then, $\left(\alpha_{\vea_n'}\hat{\ga}\ga\inv\alpha_{\vea_n'}\inv\right)\eda_1$ is an edge of $\arbmax$. Then, by the claim, we deduce that $\alpha_{\vea_n'}\hat{\ga}\ga\inv\alpha_{\vea_n'}\inv\gx$ has critical number strictly smaller than the critical number of $[1\veacero,\gx\veacero]_\arbol$, which is $s$.\\

From \Cref{bass serre para interseccion eq el elemento esta en factores 3} we have $\alpha_{\vea_n'}\hat{\ga}\ga\inv\alpha_{\vea_n'}\inv\in\gen{\gr_0,\grb_0}$.  Notice that $\alpha_{\vea_n'}\hat{\ga}\ga\inv\alpha_{\vea_n'}\inv\in \gen{\generadores}\subseteq\ce$ where the first containment is \Cref{bass serre para interseccion eq el elemento esta en gen de S 3} and the second is consequence of \Cref{bass serre para interseccion eq generadores contenidos en la interseccion}. We deduce that $\alpha_{\vea_n'}\hat{\ga}\ga\inv\alpha_{\vea_n'}\inv\in\ce\cap\gen{\gr_0,\grb_0}$. Applying that $\gx\in\ce\cap\gen{\gr_0,\grb_0}$ by hypothesis, we deduce that $\alpha_{\vea_n'}\hat{\ga}\ga\inv\alpha_{\vea_n'}\inv\gx\in\ce\cap\gen{\gr_0,\grb_0}$. Then, as we know that $\alpha_{\vea_n'}\hat{\ga}\ga\inv\alpha_{\vea_n'}\inv\gx$ has critical number strictly smaller than $s$, by induction hypothesis, we obtain that $\alpha_{\vea_n'}\hat{\ga}\ga\inv\alpha_{\vea_n'}\inv\gx\in \gen{\generadores}$. As we know that $\alpha_{\vea_n'}\hat{\ga}\ga\inv\alpha_{\vea_n'}\inv\in \gen{\generadores}$ from \Cref{bass serre para interseccion eq el elemento esta en gen de S 3}, we deduce that $\gx\in \gen{\generadores}$. This concludes the Case III.\\

\item[\textbf{Case IV}]  Suppose that the first apparition of a critical element in $[1\veacero,\gx\veacero]_\arbol$ is a critical edge, and suppose that the previous vertex is an element of the form $\gy\veagrb$ for some $\gy\in\gr*\grb$. This case is analogous to Case III.
\end{itemize}
\end{proof}
\end{teo}

\begin{obs}
Notice, on the proof of the previous result, for such $\ce$ and for such $\gr_0$ and $\grb_0$ satisfying $\fino\subseteq \gen{\gr_0,\grb_0}$, we have that the Kurosh rank of $\ce\cap\gen{\gr_0,\grb_0}$ is bounded by the Kurosh rank of $\ce$. Also notice that, on that proof, one can consider 
$$\fino=\{\zeta_\indi \tq \indi\in\fini\}\cup \{\alpha_\indj \tq \indj\in\finj\}\cup \{\beta_\indk \tq \indk\in\fink\}$$
instead of $$\fino=\{\zeta_\indi \tq \indi\in\fini\}\cup \{\alpha_\indj \tq \indj\in\finj\}\cup \{\beta_\indk \tq \indk\in\fink\}\cup \{\gx\in\gr*\grb\tq \gx\edagr\in \edges\arbmax\cup\aristasfuera^{\pm}\;\mbox{ or }\;\gx\edagrb\in \edges\arbmax\cup\aristasfuera^{\pm}\},$$ as one can prove that any  $\gr_0\subseteq \gr$ and $\grb_0\subseteq\grb$ such that $\gen{\gr_0,\grb_0}$ contains the smaller finite set then it contains the bigger one.
\end{obs}

\noindent {\bf Acknowledgements}. The author thanks Yago Antolín for the supervision and for all the crucial suggestions and comments about the paper.

\end{document}